\pdfoutput=1
\documentclass[11pt]{article}

\usepackage[margin=1.15in]{geometry}
\usepackage[T1]{fontenc}
\usepackage[utf8]{inputenc}
\usepackage{microtype}

\usepackage{mathtools}               
\usepackage{amssymb}
\usepackage{amsfonts}
\usepackage{amsopn}
\usepackage{mathrsfs}                
\usepackage{bbm}                     

\usepackage{graphicx}
\usepackage[caption=false]{subfig}   
\usepackage{array}
\usepackage{comment}                 
\usepackage[shortlabels]{enumitem}   
\setlist[enumerate]{leftmargin=.5in}
\setlist[itemize]{leftmargin=.5in}

\usepackage[title]{appendix}

\usepackage{tikz}
\usepackage{pgfplots}
\pgfplotsset{compat=1.18}
\usepgfplotslibrary{fillbetween}

\usepackage{hyperref}
\hypersetup{
  colorlinks=true,
  allcolors=[rgb]{0.10,0.25,0.55},
  pdftitle={Learning to Integrate},
  pdfauthor={O. G. Ernst, H. Gottschalk, T. Kowalewitz and P. Kruger}
}
\usepackage[all]{hypcap}
\usepackage[capitalise,nameinlink]{cleveref}

\crefname{figure}{Figure}{Figures}
\Crefname{figure}{Figure}{Figures}
\crefformat{equation}{(#2#1#3)}
\Crefformat{equation}{(#2#1#3)}
\crefrangeformat{equation}{(#3#1#4)--(#5#2#6)}
\crefmultiformat{equation}{(#2#1#3)}{ and~(#2#1#3)}{, (#2#1#3)}{, and~(#2#1#3)}

\numberwithin{equation}{section}

\newcommand{\email}[1]{\href{mailto:#1}{\nolinkurl{#1}}}
\newenvironment{keywords}%
  {\par\medskip\noindent\textbf{Key words. }\ignorespaces}{\par}
\newenvironment{MSCcodes}%
  {\par\medskip\noindent\textbf{MSC codes. }\ignorespaces}{\par}

\newcommand{\vc}{\boldsymbol c}

\newcommand{\vf}{\boldsymbol f}
\newcommand{\vg}{\boldsymbol g}

\newcommand{\vk}{\boldsymbol k}
\newcommand{\vell}{\boldsymbol \ell}
\newcommand{\vm}{\boldsymbol m}
\newcommand{\vn}{\vec{n}}
\newcommand{\vr}{\boldsymbol r}
\newcommand{\vs}{\boldsymbol s}
\newcommand{\vt}{\boldsymbol t}
\newcommand{\vu}{\boldsymbol u}
\newcommand{\vv}{\boldsymbol v}

\newcommand{\vy}{\boldsymbol y}
\newcommand{\vz}{\boldsymbol z}

\newcommand{\vA}{\boldsymbol A}
\newcommand{\vB}{\boldsymbol B}

\newcommand{\vI}{\boldsymbol I}

\newcommand{\vO}{\boldsymbol O}

\newcommand{\spS}{\mathscr{S}}

\DeclareMathOperator{\one}{{\mathbbm{1}}}

\newcommand{\valpha}{{\boldsymbol \alpha}}

\newcommand{\veta}{{\boldsymbol \eta}}
\newcommand{\vtau}{{\boldsymbol \tau}}
\newcommand{\vphi}{{\boldsymbol \phi}}

\newcommand{\vsigma}{\vec{\sigma}}

\newcommand{\vxi}{{\boldsymbol \xi}}

\newcommand{\vnu}{{\boldsymbol \nu}}

\newcommand{\ev}[1]{ {\boldsymbol{\mathsf  E}} \left[  #1 \right]}
\newcommand{\evdist}[2]{ {\boldsymbol{\mathsf  E}}_{#1} \left[  #2 \right]}

\DeclareMathOperator{\Cov}{\boldsymbol{\mathsf {Cov}}}

\newcommand{\pdist}[1]{{\mathsf #1}}

\DeclareMathOperator{\trace}{tr}
\DeclareMathOperator*{\esssup}{ess\,sup}
\DeclareMathOperator*{\essinf}{ess\,inf}

\renewcommand{\d}[1]{\ensuremath{\operatorname{d}\!{#1}}}
\newcommand{\e}{\mathrm e}
\newcommand{\R}{\mathbb{R}}

\newcommand{\oge}[1]{#1}

\title{Learning to Integrate: High-Dimensional Quadrature\\
       Using Flow-Based Generative Models}

\author{Oliver G. Ernst\thanks{Department of Mathematics, TU Chemnitz, Germany
  (\email{oernst@math.tu-chemnitz.de}, \email{toni.kowalewitz@math.tu-chemnitz.de}).}
\and Hanno Gottschalk\thanks{Institut für Mathematik, TU Berlin, Germany
  (\email{gottschalk@math.tu-berlin.de}, \email{krueger@math.tu-berlin.de}).}
\and Toni Kowalewitz\footnotemark[1]
\and Patrick Krüger\footnotemark[2]}

\date{}   

\begin{document}
\maketitle

\begin{abstract}
We address uncertainty quantification for a generic input distribution to some resource-intensive simulation, e.g., requiring the solution of a partial differential equation. 
While efficient numerical methods exist to compute integrals for high-dimensional Gaussian and other separable distributions based on sparse grids (SG), input data arising in practice often does not fall into this class. 
We therefore employ transport maps to transform complex distributions to \oge{multivariate} standard normals. 
In generative learning, a number of neural network architectures have been introduced that accomplish this task approximately. 
Examples are affine coupling flows (ACF) and ordinary differential equation-based networks such as conditional flow matching (CFM). 
To compute the expectation of a quantity of interest, we numerically integrate the composition of the inverse of the learned transport map with the simulation code output. 
As this map is integrated over a multivariate Gaussian distribution, SG techniques can be applied. Viewing the images of the SG quadrature nodes as learned quadrature nodes for a given complex distribution motivates our title. 
We demonstrate our method for monomials of total degrees for which the unmapped SG rules are exact. 
We also apply our approach to the stationary diffusion equation with coefficients modeled by exponentiated Lévy random fields, using modal expansions of Karhunen--Loève type with 9 and 25 modes.  
In a series of numerical experiments, we investigate errors due to learning accuracy, quadrature, statistical estimation, truncation of the modal series of the input random field, and training data size for three normalizing flows (ACF, conditional flow matching and optimal transport flow matching). 
We discuss the mathematical assumptions on which our approach is based and demonstrate its shortcomings when these are violated.     
\end{abstract}

\begin{keywords}
Learned sparse grid quadrature for general distributions $\bullet$
normalizing flows $\bullet$
flow matching $\bullet$
random diffusion equation $\bullet$ 
Lévy random coefficients $\bullet$ 
Smolyak sparse Gauss--Hermite  quadrature
\end{keywords}

\begin{MSCcodes}
65C30, 
65D40, 
68T07. 
\end{MSCcodes}

\section{Introduction}

Uncertainty quantification (UQ) is crucial for assessing the validity of results obtained by numerical simulation. 
Usually, any numerical simulation is based on input parameters $\veta$ entering a simulation code,  which then computes some quantity of interest (QoI), which we shall denote by $Q(\veta)$. 
Uncertainty in finite-dimensional inputs $\veta$ can be modeled by an $\R^M$-valued random variable $\veta$ with distribution $p(\veta)\d{\veta}$. 
Computing statistics of the resulting random QoI $Q(\veta)$, say, the expected value
\begin{equation} \label{eq:expectedValue}
    \evdist{\veta\sim p(\veta)}{Q(\veta)} 
    =
    \int_{\R^M} Q(\veta)\, p(\veta) \d{\veta}
\end{equation}
usually leads to problems of numerical quadrature of a potentially high-dimensional integral. 
If the density $p(\veta)$ factorizes, i.e., $p(\veta)=\prod_{j=1}^M p_j(\eta_j)$, one can approximate the integral with standard product quadrature rules obtained by tensorizing univariate rules derived from the marginal densities $p_j(y_j)$, e.g., product Gauss--Hermite quadrature \cite{Trefethen2022,KazashiEtAl2023}. 
Unfortunately, (a) most real-world distributions involve complicated dependencies and fail to factorize and (b) even if they do, the \emph{curse of dimensionality} makes the evaluation intractable in high dimensions $M\gg 1$ as the number of required function evaluations for $Q$ grows exponentially in $M$ \cite[Chapter~3]{StroudSecrest1966}. 
This is especially problematic if evaluating $Q(\veta)$ involves running an expensive numerical simulation, requiring days or weeks to complete. 

Sparse grid (SG) quadrature is a strategy to address problem (b) if the QoI $Q$ depends sufficiently smoothly on $\veta$. 
Sparse quadrature rules only require $\mathcal{O}(M\log M)$ evaluations of the integrand and therefore effectively mitigate the curse of dimensionality. 
SG quadrature rules have been constructed for numerous distributions, and in particular the standard normal distribution with density
 $\varphi(\vxi) = (2\pi)^{-M/2}\e^{-\frac{1}{2}|\vxi|^2}
                = \prod_{j=1}^M (2\pi)^{-1/2} \e^{-\frac{1}{2}\xi_j^2}$. 
However, if problem (a) remains unresolved, particularly if $p(\veta)$ is available only indirectly in the form of observational data $\mathcal{D}=(\veta_1,\ldots,\veta_N)$ sampled from $p(\veta)$, sparse grid quadrature cannot be immediately applied.

A way to tackle the problem nevertheless is to transport the problem-specific distribution $p(\veta)$ to a simpler reference distribution, e.g., the standard multivariate normal distribution \cite{santambrogio2015optimal,villani2008optimal}. 
Let thus $\vf\colon\R^M \to \R^M$ denote a bijective transport map that transports $p(\veta)$ to $\varphi(\vxi)$, i.e., the \emph{normalizing map} $\vf$ is constructed such that for $\veta\sim p$, we have $\vxi = \vf(\veta)\sim \varphi$. 
Using the change of variables formula for probability densities, this amounts to 
\begin{equation} \label{eq:transformationofDensities}
    p(\veta) = \varphi(f(\veta)) \left| \det (D\vf)(\veta) \right|,
\end{equation}
where $D\vf$ denotes the Jacobian of the map $\vf$. 
Let $\vg$ denote the inverse of $\vf$, also called the \emph{generative map} since $\vxi\sim \varphi$ implies $\vg(\vxi)\sim p$. 
In measure-theoretic terms, the generative map $\vg$ is the pushforward of the probability measure with density $\varphi$ to that associated with the density $p$.
The computation of the expected value then amounts to
\begin{equation} \label{eq:expectedValue2}
    \evdist{\veta\sim p}{Q(\veta)}
    =
    \int_{\R^M} Q(\veta) \, p(\veta) \, \d{\veta}
    =
    \int_{\R^M} Q(\oge{\vg}(\vxi)) \, \varphi(\vxi) \d{\vxi}
    \approx 
    \sum_{j=1}^{n_\text{quad}} w_j (Q \circ \vg)(\vxi_j),
\end{equation}
where on the right we have applied a generic quadrature rule with $n_\text{quad}$ weights $w_j$ and nodes $\vxi_j$ which will be specialized to a Smolyak sparse Gauss--Hermite formula in \cref{sec:SG} below.
\oge{Combining a change of variables with quadrature rules has a long history, see e.g.\ \cite{StroudSecrest1966,SagSzekeres1964,TakahasiMori1973} and, more recently, \cite{HaleTrefethen2008,GnewuchEtAl2025}.
The novel aspect of our approach lies in leveraging the potential of advanced generative machine learning to construct the mappings, thus enabling uncertainty quantification involving complex high-dimensional probability distributions which may only be accessible through samples.
}

While the theory of optimal transport guarantees existence of transport maps $\vf$ and generative maps $\vg$ \cite{santambrogio2015optimal,villani2008optimal} provided $p(\veta)>0$ on $\R^M$, there remains the problem of how to construct them, especially if $p(\veta)$ is only known through observed samples.
Fortunately, a number of solutions to this problem have recently been proposed in the context of generative learning. 
Here, approximate normalizing maps $\vf$ and generative maps $\vg$ are represented by deep neural networks that transform $p(\veta)$ layer by layer to $\varphi(\vxi)$ and are therefore called \emph{normalizing flows}. 
One prerequisite is that these networks be invertible neural networks (INN), which somewhat restricts their architectures to, e.g., \emph{affine coupling flows} (ACF) \cite{dinh2016density} or \emph{LU-nets} \cite{chan2023lu}. 
Note that both families of networks are sufficiently expressive to generate or normalize any continuous distribution \cite{teshima2020coupling,rochau2024new}. 
\emph{Continuous normalizing flows} map the target distribution $p(\veta)$ to the source distribution $\varphi(\vxi)$ by way of the flow map of a system of ordinary differential equations (ODEs). 
Here the vector field generating the exact flow is approximated by a neural network, a setting known as a \emph{neural ODE} \cite{chen2018neural}. 
Various training algorithms exist to estimate $\hat{\vf} \approx \vf$ and $\hat{\vg} = \hat{\vf}^{-1} \approx\vg$ from data, namely likelihood-based methods in the case of INNs and neural ODEs or the approximation of certain transporting vector fields like conditional flow matching (CFM) \cite{lipman2022flow} or \emph{optimal transport conditional flow matching} (OTCFM) \cite{tong2024improving}.

\begin{figure}[t]
    \centering
    \includegraphics[width=\linewidth]{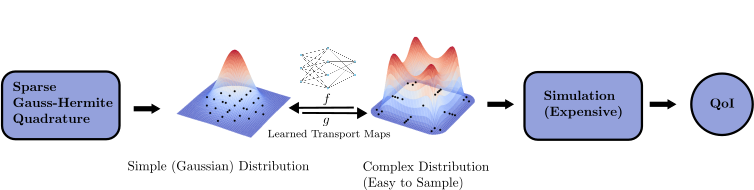}
    \caption{Overview of the learning to integrate method. The target distribution is normalized with the learned flow map $\vf$. For integration with respect to an i.i.d\ normal density, Smolyak sparse Gauss--Hermite rules allow for efficient quadrature. The associated quadrature nodes are transformed via the generative map $\vg=\vf^{-1}$ and then fed into the simulation code which computes the quantity of interest (QoI). The results, weighted with the sparse Gauss--Hermite weights, yield an approximation of the expectation of the QoI with respect to the original complex distribution.}
    \label{fig:LtI}
\end{figure}

Here we propose to learn normalizing and generative maps using cutting edge generative models from machine learning and replace the transported quadrature points 
$\veta_j := \vg(\vxi_j)$ on the right hand side of \eqref{eq:expectedValue2} by their learned counterparts $\hat{\veta}_j := \hat{\vg}(\vxi_j)$ and in this sense use learned quadratures, which explains our title \emph{learning to integrate}. 
This gives
\begin{equation} \label{eq:expectedValue3}
    \evdist{\veta\sim p(\veta)}{Q(\veta)}
    \approx 
    \sum_{j=1}^{n_\text{quad}} w_j \, Q(\hat{\vg}(\vxi_j))
    =:
    \sum_{j=1}^{n_\text{quad}} w_j \, Q(\hat{\veta}_j).
\end{equation}
Notably, the quadrature weights $w_j$ need not be adjusted. 
\Cref{fig:LtI} provides a schematic overview of our computational scheme. 
We note that our approach is particularly efficient if sampling from the data distribution $p(\veta)$ is inexpensive and $\hat \vf$ and $\hat \vg$ can be learned efficiently. 
We note that a similar approach has been taken using a learned Knothe--Rosenblatt transformation and normalizing to the uniform distribution on the unit cube in $\R^M$, see \cite{santambrogio2015optimal}. 
The Knothe--Rosenblatt map, due to its restrictive triangular structure, however, is harder to learn, and the learning algorithms we employ here have been shown to scale to dimensions in the thousands.

In this paper we implement the above strategy numerically. 
In particular, we consider the standard model problem from the uncertainty propagation literature, i.e., the stationary diffusion equation with conductivity modeled by a random field. As recently proposed \cite{ErnstEtAl2021TR}, we consider exponentiated smoothed L\'evy random field conductivities, which extend the more standard class of lognormal random fields. 
We note that multivariate Lévy distributions, due to higher-order cumulants, generally fail to tensorize if correlations are removed by a linear transformation, see, e.g., \cite{albeverio1996convoluted}. Finite-dimensional modal expansions for these distributions can be obtained similar to the Karhunen--Loève expansion \cite{ErnstEtAl2021TR}. 
For such truncated modal expansions employing 9 and 25 modes, we test the proposed algorithm on two separate test cases. 
On the one hand we integrate monomials where Smolyak sparse Gauss--Hermite quadrature gives exact results with our learned quadrature rules and report averaged errors. 
On the other hand, we compute the flux over an outflow boundary of the random diffusion equation as a QoI based on a finite element approximation of solution realizations. 
In both cases, we compare with naive MC simulations based on  $10^5$ samples. 
We observe reasonable convergence behavior whenever the underlying mathematical assumptions---most importantly a continuous distribution of the input data with support on full space---is satisfied. \oge{In this case, our proposed method significantly reduces the error and the standard deviation when compared to a Monte Carlo (MC) baseline with the same number of approximate PDE solves.}

Furthermore, we provide extensive comparative studies for the proposed algorithm. 
Besides 
 investigating the convergence of the finite element simulations and the modal expansion, we also investigate the effect of the learning sample size and the efficiency of three generative methods, ACF, CFM and OTCFM. 
 We also provide examples of failure in cases where the mathematical assumptions underlying the transport theory and sparse quadrature are violated. 
 The implementations of the data generation, machine learning and finite element approximation subtasks are publicly available\footnote{\url{https://github.com/ToniKowalewitz/LearningToIntegrate}}.

While the work presented here is focused on a proof-of-concept and extensive numerical studies,   the companion paper \cite{GottschalkEtAl2025-PP} gives rigorous proofs of consistency of our learning to integrate approach in a more restricted setting than ours, namely that the distributions have compact supports (which excludes multivariate Gaussians) and Hölder continuous densities bounded away from zero and neural ODEs are constructed with maximum-likelihood training.

The remainder of the paper is organized as follows: In \cref{sec:Gen_Models} we give a brief account of the employed generative models and in \cref{sec:SG} we recall the essentials of sparse quadrature methods. 
\Cref{sec:NumSol} explains Lévy random fields, their modal expansion, simulation and our finite element approach.  
Some technicalities on the discretization both of random fields and the PDE are given in the appendix. 
\Cref{sec:experiments} presents the numerical results and gives a systematic computational analysis of all sources of error.  
In \cref{sec:conclusions} we conclude and point out future directions of research.

\section{Generative models} \label{sec:Gen_Models}

In this section we provide the theoretical background for the neural network-based generative models to be trained to transform complex distributions (gamma/bigamma/Poisson) with not necessarily independent components into an i.i.d.\ Gaussian random vector. 
The trained models are then applied in the reverse direction to map the nodes of sparse Gauss--Hermite quadrature rules, thus yielding efficient quadrature for the respective complex distributions, where the quadrature weights remain unchanged. 
Flow matching (FM) models were, among others, introduced by \cite{lipman2022flow} and present an efficient way of training NeuralODEs (NODEs, \cite{chen2018neural}), which in turn can be interpreted as a continuation of normalizing flow (NF) models as introduced by \cite{dinh2016density}. 
This section presents an overview of normalizing flows (\cref{sec:NF}), followed by NeuralODEs (\cref{sec:cont_flows}) and  the efficient training of NeuralODEs by the flow matching approach. 

\subsection{Normalizing flows} \label{sec:NF}

In what follows we denote by $\veta = \veta_1\in\R^M$ a random variable that follows a complex probability distribution $\veta_1 \sim p_{\veta_1} = p_1$, and by $\vxi = \veta_0 \in \R^M$ a random variable that follows a known simple distribution such as a standard multivariate Gaussian $\veta_0 \sim p_{\veta_0} = p_0$.\\
 \\
In the context of generative learning, a NF usually consists of a composition of invertible, differentiable functions $\vf_i^{\theta}$ with learnable parameters $\theta$ that are trained to transform the unknown distribution $p_1$ into the known $p_0$ on a per-sample basis:
\begin{equation} \label{eq:normalizing_flow}
    \vxi = \veta_0 = \vf(\veta_1) 
    = 
    (\vf_1^{\theta}\circ\cdots\circ \vf_n^{\theta})(\veta_1).
\end{equation}
The trained NF can then be used to generate new samples from $p_1$ by sampling $\vxi \sim p_0$ from the known distribution and applying $\vg := \vf^{-1}$, giving
\begin{equation} \label{eq:normalizing_flow_reverse}
    \veta_1 
    = 
    \vf^{-1}(\veta_0) = \vg(\veta_0) 
    = 
    (\vg_n^{\theta}\circ\dots\circ \vg_1^{\theta})(\veta_0), \qquad \veta_0\sim p_0,
\end{equation}
with $\vg_i^{\theta} = (\vf_i^{\theta})^{-1}$. 
A key feature of NF is the ability of exact density estimation of $p_1$ via the change of variables formula, see \cref{eq:transformationofDensities}, which in turn allows for likelihood-based training. 
Given a data set of samples $\mathcal{D}=(\veta_1^{(i)})_{i=1}^m$ from $p_1$, one can then calculate the log likelihood of the data $\mathcal{D}$ given the parameters $\theta$ of $\vf$ as
\begin{equation} \label{eq:nf_density_est}
    \log p_1(\mathcal{D}|\theta)
    = 
    \sum_{i=1}^m \log p_1(\veta_1^{(i)}|\theta)
    =
    \sum_{i=1}^m 
    \log p_0(\vf(\veta_1^{(i)}|\theta)) + \log |\det D\vf(\veta_1^{(i)}|\theta) |.
\end{equation}
\noindent
While NF have been successfully applied in image generation (\cite{kingma2018glow},\cite{dinh2016density}) as well as inverse problems in, e.g., engineering applications \cite{krueger2025generative}, a major drawback of such models is the computational cost of calculating the Jacobian determinants of $\vf$ which scales as $\mathcal{O}(LD^3)$ where $L$ and $D$ denote the number and dimension of hidden layers of $\vf$, respectively. 
To mitigate this, as well as to achieve invertibility with tractable computational effort, in INNs based on ACF subfunctions $\vf_i$ \cref{eq:normalizing_flow_reverse} are constructed from so-called \emph{coupling blocks} based on partitioning the input $\veta \in \R^M$ into $\vk_1 := \veta_{1: i}$ and $\vk_2 := \veta_{i+1:M}$. 
Outputs $\vell_1 = \veta_{1:i}$ and $\vell_2 = \veta_{i+1:M}$ are then obtained via
\begin{equation} \label{eq:coupling_law}
    \vell_1 := \vk_1, \hspace{1cm}
    \vell_2 := \vc(\vk_2;\vm(\vell_1)).
\end{equation}
In \cref{eq:coupling_law}, $\vc$ describes an invertible \emph{coupling law} which transforms $\vk_2$ given the output of a function $\vm$ (usually a neural network) evaluated on $\vell_1 = \vk_1$. 
The input can thus be recovered by inverting the mapping, i.e.,
\begin{equation} \label{eq:add_coupling_law}
        \vk_1 = \vell_1, \hspace{1cm}
        \vk_2 = \vc^{-1}(\vell_2; \vm(\vell_1)).
\end{equation}

\noindent
The most commonly used coupling law is the affine coupling \cite{dinh2016density}
\begin{align} \label{eq:affine_coupling_law}
\begin{split}
        \vell_2 &= \vk_2\odot \exp(\vs(\vell_1)) + \vt(\vell_1),\\
        \vk_2   &= (\vell_2 - \vt(\vell_1))\odot \exp(-\vs(\vell_1)),
\end{split}
\end{align}
where $\odot$ denotes the Hadamard (componentwise) product, the exponential function is evaluated componentwise and the functions $\vs,\vt \colon \mathbb R^i \to \mathbb R^{M-i}$
effect a scaling and translation of $\vell_1$ in \cref{eq:affine_coupling_law}.\\
\ \\
Coupling blocks result in triangular Jacobian matrices, which reduces computational cost in, e.g., evaluating \cref{eq:nf_density_est}. 
Furthermore, as coupling blocks are invertible by construction, the functions $\vs,\vt$ comprising $\vm$ can be chosen arbitrarily, for instance as simple, not necessarily invertible feedforward neural networks. 
While INNs based on ACF possess universal approximation properties \cite{teshima2020coupling} for bijective mappings under rather mild conditions, they may in practice still suffer from their architectural constraints.

\subsection{Continuous flows: neural ODE and conditional flow matching} \label{sec:cont_flows}

Residual neural networks (ResNets) \cite{HeEtAl2016} update hidden states by repeatedly adding learned corrections,
\[
\vy_{t+1} = \vy_t + \vr_t(\vy_t,\theta_t),
\]
which can be interpreted as a forward Euler discretization of a continuous-time dynamical system. Taking the step size to zero leads to a \emph{neural ordinary differential equation (neural ODE)} \cite{chen2018neural},
\[
\frac{\d \vy(t)}{\d t} = \vv(\vy(t),t;\theta),
\]
where the hidden state \(\vy(t)\) evolves according to a time-dependent vector field \(\vv\) represented by a neural network with parameters \(\theta\). Given an initial state \(\vy(0)=\vy_0\), the output of the Neural ODE is defined as \(\vy(1)\), obtained by integrating the ODE from \(t=0\) to \(t=1\). In practice, this integration is carried out using a numerical ODE solver, which adaptively evaluates the vector field \(\vv\) along the trajectory.

When Neural ODEs are used to transform random variables, the associated flow map \(\vphi_t\) defines a continuous and invertible transformation of the input space. Applying this flow to an initial random variable \(\veta_0\sim p_0\) induces a continuous family of probability densities \(\{p_t\}_{t\in[0,1]}\) via pushforward,
\[
p_t(\veta)
=
([\vphi_t]_* p_0)(\veta)
=
p_0\!\left(\vphi_t^{-1}(\veta)\right)
\left|
\det D_\veta \vphi_t^{-1}(\veta)
\right|.
\]
Under mild regularity assumptions on the vector field \(\vv\), the evolution of the density along the flow satisfies the instantaneous change-of-variables formula \cite{chen2018neural},
\begin{equation} \label{eq:cont_cow}
    \frac{\partial\log p_t(\veta)}{\partial t}
    =
    -\trace(D_\veta \vv(\veta,t)\mid_{\veta = \veta(t)})
\end{equation}
This identity forms the basis of \emph{continuous normalizing flows} (CNFs), which are typically trained by maximizing likelihoods computed through numerical integration of both the state and the log-density. While this framework allows for flexible and expressive models, the reliance on ODE solvers during training leads to substantial computational overhead.\\
\ \\
\noindent
Flow matching (FM) was introduced as a simulation-free alternative for learning such continuous probability flows \cite{lipman2022flow,onkenRuthotto2020-PP,liu2022flow,albergo2022building}. Instead of learning the vector field indirectly via likelihood maximization and repeated numerical integration, FM directly regresses a neural vector field \(\vv_t\) toward a target vector field \(\vu_t\) that generates a prescribed probability path \(p_t\). The flow matching objective is
\[
    \mathcal{L}_{\mathrm{FM}}(\theta)
    =
    \evdist{t\sim\pdist{U}(0,1),\,\veta\sim p_t}
           {\|\vv_t(\veta;\theta)-\vu_t(\veta)\|^2}.
\]
\noindent
In general, the target vector field \(\vu_t\) is intractable to compute explicitly. To address this, the probability path \(p_t\) is represented as a mixture of conditional paths,
\[
p_t(\veta)=\int p_t(\veta\mid\vz)\,q(\vz)\,\d\vz,
\]
where \(\vz\) is a conditioning variable drawn from a known distribution \(q\). If each conditional density \(p_t(\veta\mid\vz)\) is generated by a corresponding conditional vector field \(\vu_t(\veta\mid\vz)\), then the marginal vector field generating \(p_t\) satisfies
\[
    \vu_t(\veta)
    =
    \evdist{\vz\sim q}{\frac{\vu_t(\veta\mid\vz)\,p_t(\veta\mid\vz)}{p_t(\veta)}}.
\]
Although this expression remains intractable, it can be shown that minimizing the flow matching loss is equivalent, up to a constant independent of \(\theta\), to minimizing the \emph{conditional flow matching} (CFM) objective \cite{lipman2022flow,tong2024improving}, which regresses $\vv_t$ against the conditional vector field $\vu_t(\veta|\vz)$ on a per-sample basis: 
\begin{equation}\label{eq:loss_cfm}
    \mathcal{L}_{\mathrm{CFM}}(\theta)
    =
    \evdist{t\sim\pdist{U}(0,1),\vz\sim q,\,\veta\sim p_t(\cdot\mid\vz)}
           {\|\vv_t(\veta;\theta)-\vu_t(\veta\mid\vz)\|^2}.
\end{equation}
\noindent
A widely used choice for the conditional probability paths is a Gaussian interpolation between samples \(\veta_0\sim q_0\) and \(\veta_1\sim q_1\),
\[
    p_t(\veta\mid\veta_0,\veta_1)
    =
    \pdist{N}\!\left(
    \veta\,\middle|\,
    (1-t)\veta_0+t\veta_1,\,
    \sigma^2\vI
    \right),
\]
with conditioning variable \(\vz=(\veta_0,\veta_1)\).
For this choice, the corresponding conditional vector field admits a simple closed-form expression \cite{tong2024improving},
\[
    \vu_t(\veta\mid\vz)=\veta_1-\veta_0.
\]
Sampling \(\veta_0\) and \(\veta_1\) independently yields \emph{independent conditional flow matching}, while choosing the conditioning distribution as a 2-Wasserstein optimal transport coupling between \(q_0\) and \(q_1\) leads to optimal transport conditional flow matching (OTCFM), which potentially improves training efficiency and inference time \cite{tong2024improving}.

\section{Sparse grid quadrature} \label{sec:SG}

In this section we briefly recall Smolyak \emph{sparse grid quadrature} based on Gauss--Hermite (GH) formulas for approximating the integral \cref{eq:expectedValue3}.
In view of the discussion preceding \cref{eq:expectedValue2}, the expectation of the QoI can be approximated by applying a quadrature rule to the integrand $Q_{\vg} := Q \circ \vg$, where $\vg$ is the generative map that transports an i.i.d.\ multivariate Gaussian distribution to the multivariate Lévy distribution under consideration.
We are thus led to approximating the $M$-variate integral
\begin{equation} \label{eq:Int}
   I = \int_{\mathbb R^M} Q_{\vg}(\vxi) \varphi(\vxi) \d{\vxi}, 
\end{equation}
in terms of the standard Gaussian density $\varphi(\vxi) = (2\pi)^{-M/2} \exp(-\|\vxi\|^2/2)$.
For smooth integrands $Q_{\vg}$ and small values of $M$, product Gauss quadrature rules provide approximations to multivariate integrals \cref{eq:Int} which converge exponentially fast with increasing number of evaluations of the integrand \oge{(cf.\ \cite{KuoEtAl2017,GnewuchEtAl2025})}.
Since the number \oge{$k^M$} of function evaluations required by the product $k$-point rule in $M$ variables grows exponentially with $k$, such product rules suffer from the curse of dimensionality for large values of $M$, which in the examples to follow will take the values $M=9$ and $M=25$, for which product rules are no longer feasible.

\emph{Sparse grid} quadrature for multivariate integrals is a clever way of combining product rules each of which employs a large number of quadrature nodes in only a small number of dimensions.
Their building blocks consist of a sequence $(I_k)_{k \in \mathbb N}$ of univariate quadrature rules
\[
    I_k(q) = \sum_{j=1}^{n_k} w_j^{(k)} q(\xi_j^{(k)})
\]
with positive weights $\{w_j^{(k)}\}_{j=1}^{n_k}$ and nodes $\{\xi_j^{(k)}\}_{j=1}^{n_k}$. 
For our purposes we choose the sequence of univariate GH rules with $n_k = k$, i.e., the nodes of $I_k$ are the zeros of the (probabilist's) Hermite orthogonal polynomial of degree $k$ and the weights chosen such that $I_k$ is interpolatory (cf.\ \cite[p.\ 21]{Gautschi2004}). 
In the multivariate setting, rather than using the same number of nodes in each variable, consider product rules of the form
\begin{equation} \label{eq.quadrature-multi-index}
    I_{\vnu}(u) 
    =
    (I_{\nu_1} \otimes \dots \otimes I_{\nu_M})(u)
    :=
    \sum_{\valpha \le \vnu} w_\valpha q(\vxi_\valpha),
    \qquad\vnu = (\nu_1,\dots,\nu_M) \in \mathbb N^M,
\end{equation}
determined by a fixed multi-index $\vnu \in \mathbb N^M$, quadrature weights $w_{\valpha} = \prod_{m=1}^M w_{\alpha_m}$ for all multi-indices $\valpha$ with $\valpha \le \vnu$, where the relation is taken componentwise, and $\vxi_\valpha \in \Xi^{(\nu_1)} \times \dots \times \Xi^{(\nu_M)}$, i.e., $1 \le \alpha_m \le \nu_m$ for $m=1,\dots,M$.
The \emph{Smolyak sparse grid} construction is based on the quadrature \emph{difference} or \emph{detail} functionals
\[
    \Delta_\ell := I_\ell - I_{\ell-1},
    \quad 
    \ell \in \mathbb N, \qquad I_0 := 0,
\]
which are combined to form the \emph{Smolyak sparse grid quadrature rule} of level $L\in \mathbb N$
\[
    S_{L}(q)
    =
    S_{L,M}(q)
    :=
    \sum_{|\valpha|\le L} (\Delta_{\alpha_1} \otimes \dots \otimes \Delta_{\alpha_M})(q).
\]
An equivalent expression as a linear combination of product rules of different orders such as \cref{eq.quadrature-multi-index} was given in \cite{WasilkowskiWozniakowski1995} as
\[
    S_{L,M}(q)
    =
    \sum_{\vnu \in {\cal P}_{L,M}} 
    (-1)^{L-|\vnu|} \binom{M-1}{L-|\vnu|} I_\vnu(q),
    \qquad
    {\cal P}_{L,M} := \left\{ \vnu \in \mathbb N^M : L - M < |\vnu| \le L \right\}.
\]
\oge{
Results of Kazashi et al.\ \cite{KazashiEtAl2023} on the convergence rate of GH quadrature  have recently been extended to a convergence analysis of Smolyak sparse GH quadrature in \cite{KazashiEtAl2025}, where it is shown that, for integrands $q$ in Gaussian Sobolev space (also known as Hermite space of finite smoothness \cite{GnewuchEtAl2025}) with bounded mixed derivatives of order $r$
\begin{equation}
\label{eq:SG_Error}
    |I(q) - S_{L,M}(q)|
    \le
    \log(N)^{(M-1)(1 + \frac{r}{2})} N^{-r/2},
\end{equation}
where $N$ denotes the number of quadrature nodes occurring in $S_{L,M}$.
As pointed out in \cite{KazashiEtAl2025}, sparse GH quadrature fails to achieve the optimal rate of $N^{-r}$ achieved by certain QMC methods combined with a suitable change of variables.
We have nonetheless chosen sparse GH quadrature for our numerical experiments in \Cref{sec:experiments} for ease of implementation and because their multivariate polynomial exactness properties are used there to measure the quality of the learned transport map.
Our methodology can equally well be combined with QMC techniques which, given the suboptimality of sparse GH quadrature, will become increasingly competitive as the dimension grows.
An approach combining carefully designed transport maps with randomized QMC can be found in \cite{Liu2024}.
}

For the specific variant of sparse grid quadrature rule $S_{L,M}$ employed in the numerical experiments in \cref{sec:experiments}, i.e., the Smolyak construction based on the univariate sequence $(I_k)_{k \in \mathbb N}$ of GH rules consisting of $n_k = k$ nodes, it was shown in \cite{heiss2008likelihood} that these integrate polynomials in $M$ variables up to a total degree of $2(L-M+1)-1$.
In order for the level of sparse grid quadrature formulas to begin with 1 for a fixed dimension $M$, we shall denote the sparse grid level for fixed $M$ by the $M$-dependent index 
\[
    L_\text{SG} = L_\text{SG}(M,L) = L - M + 1, \qquad L_\text{SG} = 1,\dots,4,
\]
resulting in $L = L_\text{SG} -1 + M$ and exactness of $S_{M,L}$ for total degree up to $2L_\text{SG}-1$.

\section{Random PDEs with smoothed Lévy field coefficients} \label{sec:NumSol}

 In this section we present an extended example where non-independent non-Gaussian random input variables arise in the solution of random PDEs.
 The well-worn model problem for UQ for PDEs with uncertain inputs is the stationary diffusion equation with a random diffusion coefficient function \cite{BabuskaEtAl2010,ErnstEtAl2018,DungEtAl2023}.
 In the most common setting, the input is a Gaussian random field transformed to ensure positivity. 
 In \cite{ErnstEtAl2021TR} it was shown how to generalize the probability law of the input field to Lévy random fields smoothed by convolution with Matérn kernels to ensure continuous realizations.
 The resulting coefficient fields are represented by modal expansions with non-Gaussian non-independent random coefficients which follow a multivariate Lévy distribution.
 To provide the necessary background, we first discuss Lévy noise fields and noise fields smoothed by pseudodifferential operators following \cite{albeverio1996convoluted,ErnstEtAl2021TR}. These are then introduced into the diffusion equation as coefficients. 
We also discuss the numerical simulation of random fields and the random diffusion PDE. Some technicalities are referred to the appendix. 

\subsection{Random stationary diffusion equation}

Given a bounded and connected domain $D \subset \mathbb{R}^d$ with Lipschitz boundary $\partial D$, a measurable partition of its boundary $\partial D = \partial_D \cup \partial_N$, denoting subsets where Dirichlet and Neumann boundary conditions are imposed, respectively, such that $\partial_D \cap \partial_N = \emptyset$ and such that $\partial_D$ has positive surface measure,
we consider the boundary value problem for the stationary diffusion equation
\begin{equation} \label{eq:diff_eq}
	\left\{ 
	\begin{array}{rll}
		-\nabla \cdot (a \nabla u) = f_\text{source} & \mbox{ in }D,\\
		u = g_D                        & \mbox{ along }\partial_D,\\
		\vn \cdot a\nabla u = g_N    & \mbox{ along }\partial_N,
	\end{array}
	\right.
\end{equation}
with given coefficient function  $a \in L^\infty(D)$, 
source term $f_\text{source}\in L^2(D)$, 
Dirichlet boundary data $g_D \in H^{\frac{1}{2}}(\partial_D)$,
Neumann boundary data $g_N \in H^{-\frac{1}{2}}(\partial_N)$, 
and $\vn$ denoting the outward unit normal vector along $\partial D$. 
As usual, we interpret \cref{eq:diff_eq} in the weak sense.

The boundary value problem \cref{eq:diff_eq} models a great variety of phenomena in the physical sciences, among these single-phase saturated groundwater flow in a porous medium governed by Darcy's law, which expresses the (pointwise) volumetric flux  as a function of the hydraulic head $u$ by $-a(x)\nabla u(x)$.
In such a setting, the precise value of the conductivity coefficient is typically uncertain, e.g.\ derived from sparse information based on limited observations.
Modeling such uncertainty by introducing a probability distribution on the set of admissible coefficient functions $a$ results in a random PDE.

To ensure pathwise existence and uniqueness of solutions for problem \cref{eq:diff_eq}, we need the differential operator in the boundary value problem to be strictly elliptic for every realization $a(\cdot, \omega)$. 
The assumption 
\begin{equation} \label{ellipticity inequality}
	0 
	< 
	\essinf_{x\in D} a(x,\omega) 
	\leq 
	\esssup_{x\in D} a(x,\omega) 
	< 
	\infty
    \quad \text{ for almost all } \omega\in\Omega
\end{equation}
suffices to guarantee the existence of a unique $u=u(\cdot,\omega)\in H^1(D)$ that solves \cref{eq:diff_eq} with $a=a(\cdot,\omega)$. 

Gaussian random fields $Z(\cdot,\omega)$ violate this assumption since they take values on all of $\R$. 
Here uniform ellipticity can be achieved by applying a suitable transformation $a = T(Z)$ with uniformly positive values, e.g., $T = \text{constant} + \exp(\cdot)$, resulting in a lognormal random field.
The same transformation can be used for random fields following a bigamma distribution. 
Random fields with strictly positive realizations such as the gamma distribution can be used without such a transformation. 
A detailed discussion of existence of solutions for Lévy random fields can be found in \cite{ErnstEtAl2021TR}. 

\begin{figure}[ht]
    \centering
    \includegraphics[height=0.19\textheight]{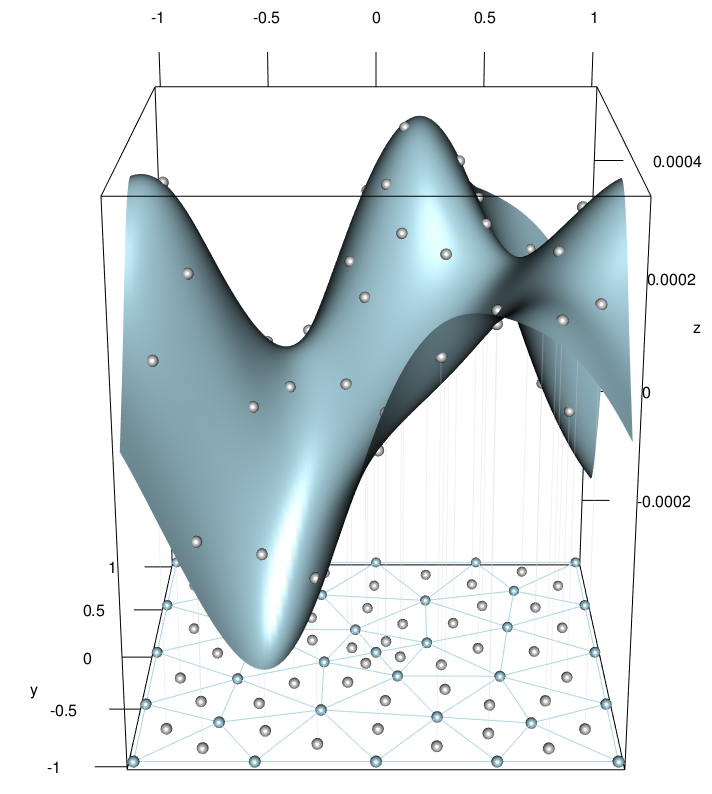}\hspace{.1cm}
    \includegraphics[height=0.19\textheight]{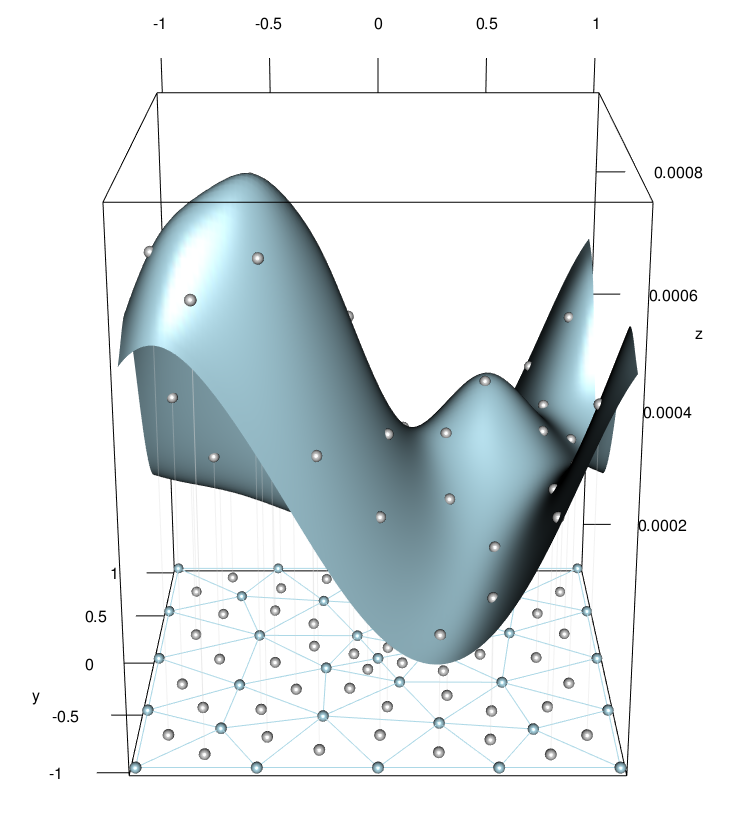} \hspace{.1cm}    
    \includegraphics[height=0.19\textheight]{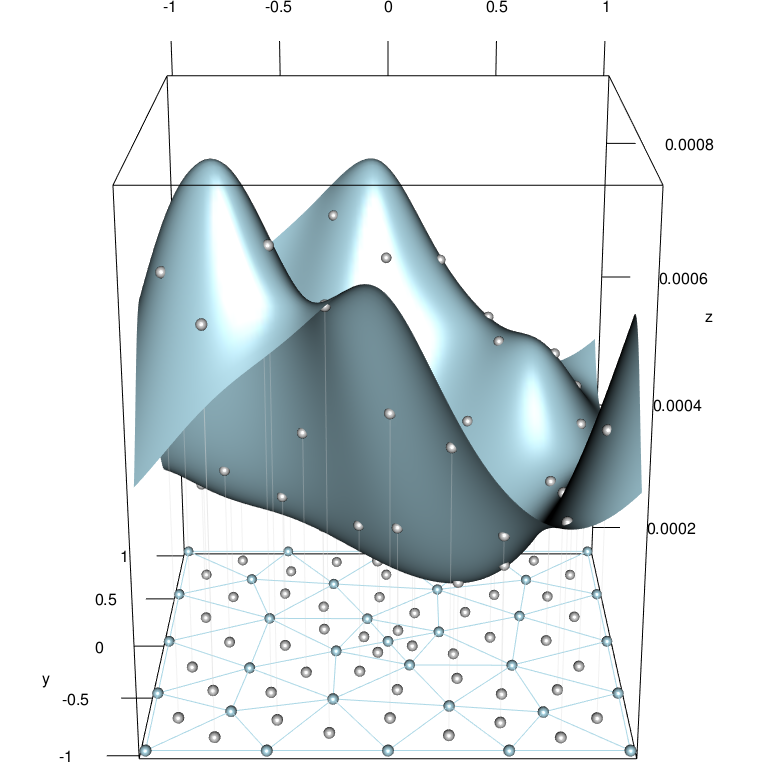}
    \hspace{.1cm}\includegraphics[height=0.19\textheight]{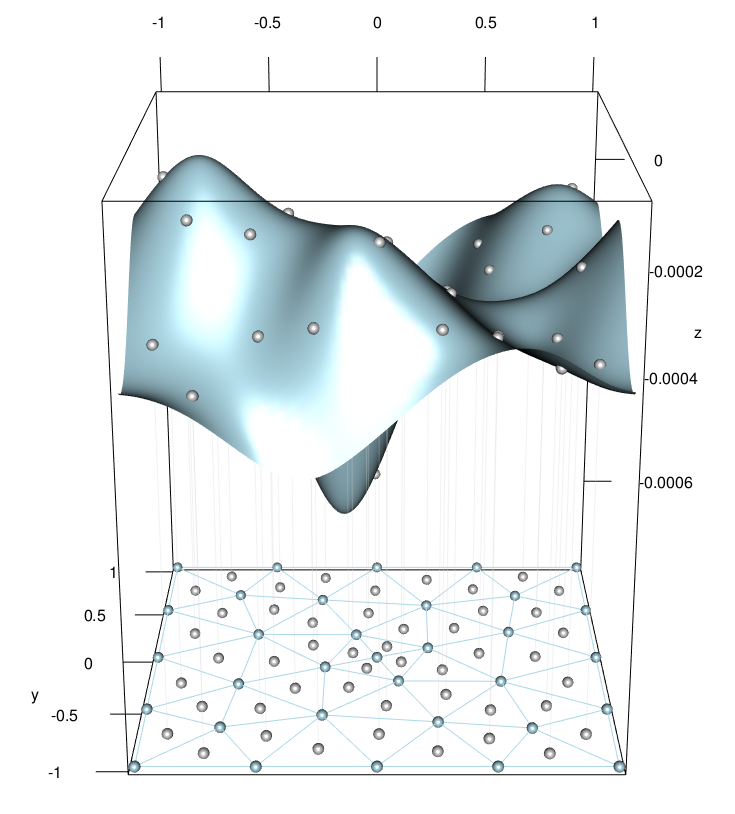}
    \caption{Realizations of smoothed Lévy noise fields with different distributions: from the left Gaussian, Poisson, gamma and bigamma interpolated to the coarse FE mesh vertices (green points) and FE quadrature points (grey points).}
    \label{fig:LevyFields}
\end{figure}

We aim to compute the expected value of a QoI of random solutions. 
Specifically, we choose the flow through a fixed part of the boundary $\partial_\text{out} \subset \partial_D$, i.e.
\begin{equation} \label{eq:QoI}
  Q(\omega) = \int_{\partial_\text{out}} -a(x,\omega)\nabla u(x,\omega)\cdot \vn\,\d s.
\end{equation}
We employ a mixed finite element discretization to approximate solutions of realizations of problem  \cref{eq:diff_eq}, since this gives better numerical control of the discretization error of the flux QoI \cref{eq:QoI}. 
The quadrature nodes where the finite element assembly procedure requires evaluation of coefficient field realizations are, in general, not aligned with the grid points $\Gamma$. 
We therefore use bilinear interpolation to interpolate from the torus grid to these quadrature nodes. 
\Cref{fig:LevyFields} shows realizations for each of the four smoothed noise fields considered.    
 For details on the FE discretization, we refer to \cref{app:MixedFE}, where also 
 the discrete approximation of the QoI is given in \eqref{eq:QoI-flux}. 

 \subsection{Lévy noise fields}
 
 Random noise fields cannot be represented as functions with values defined pointwise, 
 but can be treated using the more comprehensive theory of generalized random fields, i.e., stochastic processes indexed by a nuclear locally convex topological vector space such as the Schwartz space $\spS$ of rapidly decaying smooth functions on $\R^d$ .
 Generalized random fields may exhibit an absence of spatial correlation, as is the case with Gaussian white noise.
 
 A generalized random field on the nuclear space $\spS$ can be characterized by its \emph{characteristic functional} $\varphi: \spS \to \mathbb C$, which is continuous, positive definite and satisfies $\varphi(0) = 1$.
 Indeed, the Bochner-Minlos theorem \cite[Theorem~2.4.3]{ito1984foundations} associates with a characteristic functional $\varphi$ a random variable $Z$ taking values in the topological dual $\spS'$ of $\spS$ with probability measure $\mu$ defined on the Borel $\sigma$-algebra $\mathfrak B(\spS')$ and uniquely determined by its Fourier transform
 \[
     \chi(h) = \int_{\spS'} \e^{i Z(h)} \, \mu(\d Z), \qquad h \in \spS.
 \]
Specifically, we consider Lévy noise fields determined by their characteristic functional
\begin{equation} \label{eq:LevyNoise}
    \chi(h) = \e^{\int_{\R^d} (\psi\circ h)(x) \d{x}}, \qquad h \in \spS.    
\end{equation}
Here $\psi(t)$ is a Lévy characteristic \cite{applebaum2009levy} taking the form
\[
    \psi(t)
    =
    ibt - \frac{1}{2}\sigma^2t^2
    + \int_{\R\setminus\{0\}} \left(\e^{ist}-1-ist \one_{\{|s|\leq 1\}}\right) \nu(\d{s}),
\]
where $\nu$ is a L\'evy measure that satisfies $\int_{\R\setminus\{0\}}\min\{1,s^2\} \nu(\d{s})<\infty$, $b\in\R$ and $\sigma^2\geq 0$. 
The three parameters $(b,\sigma^2,\nu)$ are known as the \emph{characteristic triplet}. 
Under mild conditions on the existence of first moments, it has been shown \cite{ErnstEtAl2021TR} that \cref{eq:LevyNoise} are the only generalized random fields which are stationary, can be extended to characteristic functions $h=\one_\Lambda$, $\Lambda\subseteq \R^d$ bounded and measurable and have independent increments \cite{kallenberg1983random}, i.e., $Z(\one_{\Lambda_j})$ are independent whenever $\Lambda_j$ are disjoint subsets of $\R^d$.

For our numerical experiments, we utilize the characteristic functions defined by specific Lévy distributions, namely 
the Gaussian distribution as a well-studied reference \cref{eq:LevyLaws_gauss}, 
the Poisson distribution \cref{eq:LevyLaws_pois}, 
the one-sided gamma \cref{eq:LevyLaws_gamma} and 
the bigamma distribution \cref{eq:LevyLaws_bigamma}
\begin{subequations}
\begin{align} \label{eq:LevyLaws_gauss}
    \psi(t) &= -\frac{1}{2}\sigma^2 t^2,        &&\text{triplet }(0,\sigma^2,0),\\[.4em]
    \psi(t) &= \lambda \left(\e^{i t }-1\right),&&\text{triplet }(\lambda,0,\lambda\delta_1), 
               \label{eq:LevyLaws_pois}\\
   \psi(t) &= \int_{(0,\infty)} \left(\e^{i t s}-1\right) 
              \lambda \frac{\e^{-\beta s}}{s}\,\d{s}, 
    &&\text{triplet }\left(\lambda (1-\e^{-\beta}),0,\lambda\frac{\mathbbm{1}_{(0,\infty)}}{s} \d{s} \right),
               \label{eq:LevyLaws_gamma}\\
    \psi(t)
    &=
    \int_{\mathbb{R}\setminus\{0\}} 
    \left(\e^{i t s}-1\right) \lambda \frac{\e^{-\beta |s|}}{2|s|}\,\d{s},
    &&\text{triplet}\left(0,0,\lambda\frac{1}{|s|}\d s\right) \label{eq:LevyLaws_bigamma}.
\end{align}
\end{subequations}
Here $\lambda,\sigma^2, \beta >0$ and $\delta_1$ is the Dirac measure of unit mass in $s=1$.
The three parameters $\sigma^2$, $\lambda$, and $\beta$ are adjusted in order to scale the noise variance to the same values. 
To this end, we chose $\beta=1$ and set $\sigma=\sqrt{\lambda}$. 

\subsection{Smoothed Lévy fields from partial pseudodifferential equations}

Let $\Delta=\sum_{j=1}^d\partial^2_{j}$ denote the Laplacian on $\R^d$, $m^2>0$ and $\alpha>0$. 
We consider the stochastic partial pseudodifferential equation
\begin{equation} \label{eq:SPPDE}
    (-\Delta+m^2)^\alpha Z_{K_\alpha} = Z,    
\end{equation}
where $Z$ is a Lévy noise field. 
As $Z(\omega)\in\spS'$ almost surely due to the Bochner--Minlos theorem, the above equation is solved by $Z_{K_\alpha}(\omega)=(-\Delta+m^2)^{-\alpha} Z(\omega)$, since $(-\Delta+m^2)^{\alpha}$ is continuously invertible on $\spS'$ and $Z_{K_\alpha}$ is again a generalized random field. 
If the smoothing parameter $\alpha>0$ is sufficiently large, $\alpha > d+\max\left\{0,\frac{3d-12}{8}\right\}$, then the paths of $Z_{K_\alpha}$ can be realized as continuous functions, see \cite[Theorem 2.11]{ErnstEtAl2021TR}.

One way of representing $Z_{K_\alpha}$ is via the Green's function of the pseudodifferential operator $(-\Delta+m^2)^\alpha$. 
We define
\[
    K_\alpha(x)
    =
    \frac{1}{(2\pi)^d}
    \int_{\R^d}
    \frac{\e^{i\kappa\cdot x}}{(|\kappa|^2+m^2)^\alpha} \d{\kappa}
\]
and obtain that $Z_{K_\alpha}(x)=Z(K_\alpha(x-\cdot)) = (K_\alpha * Z)(x)$, where $*$ denotes convolution. 
Hence the Green's function is a Matérn function with parameters $\alpha$ and $m^2$. 
It is easily seen that the covariance function of $Z_{K_\alpha}$ is given by
\[
    \Cov[Z_{K_\alpha}(x),Z_{K_\alpha}(x')] = - \psi''(0)K_{2\alpha}(x-x'),
\]
where $\psi''(0)$ is the second derivative of the Lévy characteristic evaluated at zero, see \cite{albeverio1996convoluted,ErnstEtAl2021TR}. 
We therefore see that we can obtain statistically different random fields with the same covariance structure, see \cref{fig:LevyFields} for realizations of the random field models \crefrange{eq:LevyLaws_pois}{eq:LevyLaws_bigamma}.  

\subsection{Circulant embedding and generalized Karhunen--Lo\`eve expansion}

The simulation and expansion of Lévy random fields is based on the Mercer expansion of the smoothing function as described in \cite{bachmayr2018representations,ErnstEtAl2021TR,Reese2021}. 
As the smoothing function acts via convolution, the spectrum of the corresponding integral operator on $L^2(\R^d,\d{x})$ is continuous. 
It is therefore necessary to restrict the computational domain to a compact set in order to obtain an expansion with discrete eigenvalues. 
In this work we use the technique of circulant embedding for the compactification of the computational domain, i.e., we replace $\R^d$ with a sufficiently large torus $\mathbb{T}^d$ and replace the Laplacian $\Delta$ on $\R^d$ with its counterpart with periodic boundary conditions, which we again denote by $\Delta$. 
Likewise, we replace the Matérn kernels for smoothing on $\R^d$ with corresponding integral kernels of the pseudodifferential operators $(-\Delta+m^2)^{-\alpha}$ on $\mathbb{T}^d$. 
From the techniques used in \cite{bachmayr2018representations,ErnstEtAl2021TR,Reese2021} we can see how the error due to the circulant embedding can be efficiently controlled, and \cite{ErnstEtAl2021TR,Reese2021} show how this can be used to control the error in the $H^1$-norm of the solution of a PDE with smoothed Lévy random field coefficients.

\oge{On the torus, the pseudodifferential operator associated with the Matérn kernels has the trigonometric functions $\sin(\kappa\cdot x)$ and $\cos(\kappa\cdot x)$ as eigenfunctions with eigenvalues $(|\kappa|+m^2)^{-\alpha}$, where $\kappa\in \pi \mathbb{Z}^d$. 
By restriction to a cube of discrete nodes centered at the origin with side length $M=(2r+1)^d$, $r\in \mathbb{N}_0$, we obtain the generalized Karhunen--Loève expansion truncated at level $r$}
\begin{equation} \label{eqa:MercerMatern-indep}
    Z_{K_{\alpha,r}^c}(x)
    =
    \frac{Z(1)}{2^{d}m^{2\alpha}}
    + \frac{1}{2^{d-1}}
      \sum_{\kappa\in \pi \mathbb{Z}_{r,+}^d} 
      \frac{\cos(\kappa\cdot x)Z(\cos(\kappa\cdot (\cdot)))
             +\sin(\kappa\cdot x)Z(\sin(\kappa\cdot (\cdot)))}
           {(|\kappa|^2+m^2)^{\alpha} },
\end{equation}
\oge{
where the sum runs over multi-indices in the set $\mathbb{Z}_{r,+}$ of nonzero modes in the $d$-dimensional cube mentioned above excluding the zero mode and with the first entry to $\kappa$ positive, see  \cref{app:DiscretizationLeyField} for details}. This truncated expansion contains the vector of random modal coefficients 
\[
    \veta 
    = 
    [Z(1),
     \{Z(\cos(\kappa\cdot(\cdot)),   
       Z(\sin(\kappa\cdot(\cdot))\}_{\kappa\in\mathbb{Z}_{r,+}^d}]
    \in \R^M.
\]
It is easily seen that the distribution of the random coefficient vector $\veta$ is also of Lévy type by, e.g., observing that its characteristic function is
\begin{equation} \label{eq:LevyXi} 
    \ev{\e^{i\vtau\cdot \veta}}
    =
    \exp\left(\int_{\mathbb{T}^d} 
    \psi\left(\sum_{\ell=1}^{M}\tau_\ell \zeta_\ell(x) \right)\, \d{x}\right),
    \qquad \vtau\in \R^{M},
\end{equation}
where $\zeta_\ell$ denote the independent modes $M$ $1$, $\cos(\kappa\cdot(\cdot))$ and $\sin(\kappa \cdot (\cdot))$, $\kappa\in\pi\mathbb{Z}_{r,+}^d$, and $\psi$ is the Lévy characteristic defining $Z$, cf.\ \cref{eq:LevyNoise}. 
As we can multiply $\psi$ in this formula by an arbitrary positive constant $t>0$, we obtain characteristic functions as $l$-th roots choosing $t=\frac{1}{l}$, which is equivalent to infinite divisibility, i.e., the  distribution of $\veta$ being a Lévy distribution.  

For the discretization of the noise and efficient sampling of $\veta$ based on the fast Fourier transform of the random field $Z_{K^c_\alpha}$ on a periodic, $d$-dimensional lattice $\Gamma \subset \mathbb T^d$, we again refer to \cref{app:DiscretizationLeyField}.

\section{Numerical experiments} \label{sec:experiments}

We present the results of the complete workflow outlined in the previous sections of learning a transport map and approximating the expectation of a QoI by sparse grid quadrature. 
We first provide details of the training of our generative models, then present our results on the quadrature of monomials and the random diffusion equation with Lévy random field coefficient. 
We then systematically discuss the contribution of different sources of error and provide ablation studies.

\subsection{Training of generative models}

\begin{figure}[!ht]
    \centering
    \includegraphics[width=0.95\linewidth]{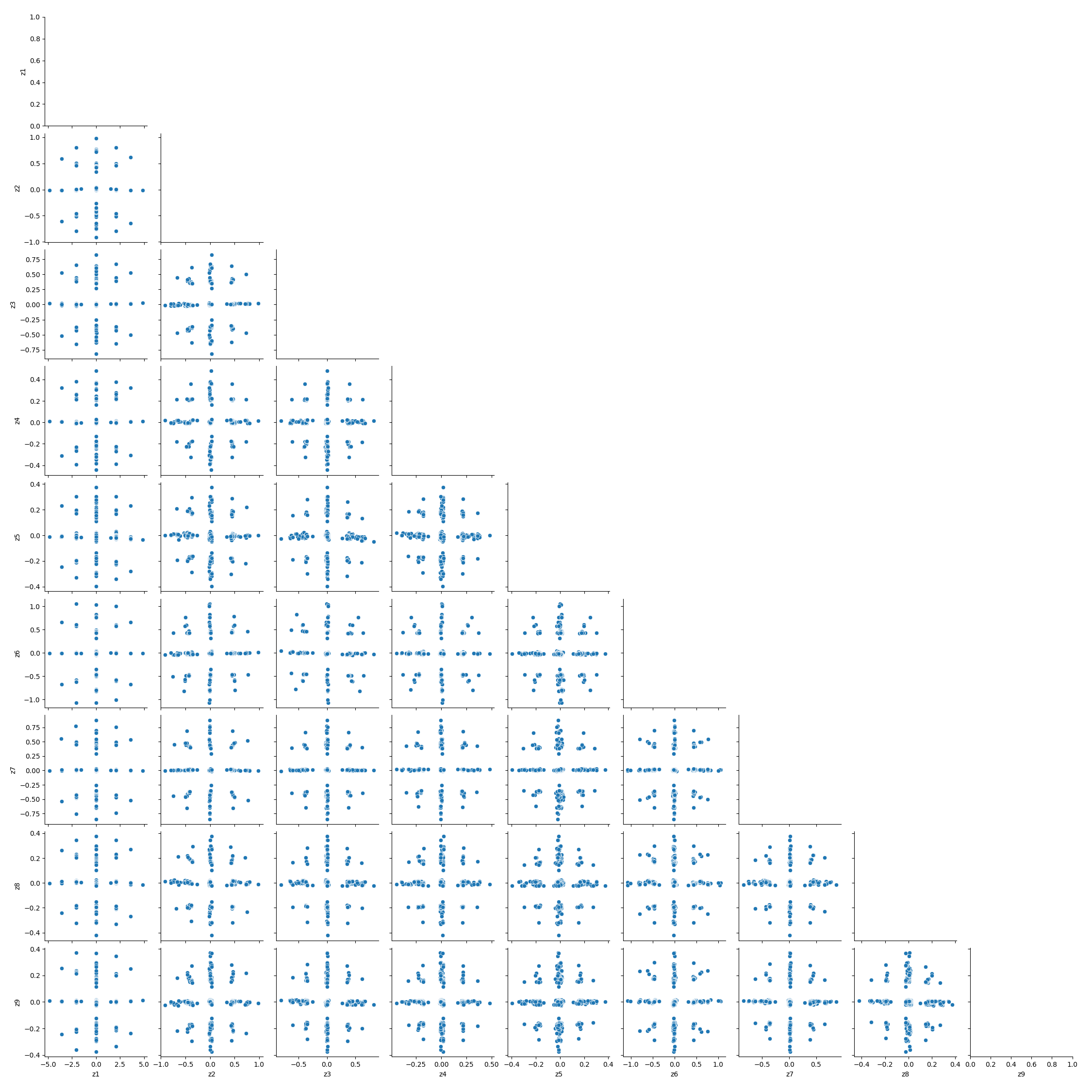}
    \caption{Images of Smolyak sparse Gauss--Hermite nodes under the transport map for the 9-mode bigamma distribution generated by the ACF model displayed in all pairs of distinct coordinate directions. 
    }
    \label{fig:Sparse_Grid_INN}
\end{figure}

We trained both ACF (\cref{sec:NF}) as well as CFM and OTCFM models (\cref{sec:cont_flows}) to transform samples from standard Gaussian distributions into the target distributions. 
To implement ACF models we used the \texttt{Framework for Easily Invertible Architectures} \cite{freia}. 
Models consisted of 4 affine coupling blocks that use the coupling law described in \cref{eq:affine_coupling_law}. 
For the sub-functions $s,t$ we chose fully connected feed-forward neural networks, each consisting of three layers with 200 neurons and Softplus \cite{zheng2015improving} activation functions in between the layers. 
Models were trained for 3000 epochs with an initial learning rate of 0.001, which was reduced by a factor of 10 after 1000 and 2000 epochs. 
Training was performed in both forward and reverse direction, i.e., transforming target distribution samples into Gaussian ones and vice versa, as this was found to improve performance. 
In both directions, Maximum Mean Discrepancy \cite{gretton2006kernel}, a kernel-based loss function that estimates the similarity of two distributions only available through samples, was applied. \\
\ \\
OTCFM models were implemented using the \texttt{torchcfm} library \cite{tong2024improving}. 
The models used to approximate the vector field $\vv_t(\veta,\theta)$ (see \eqref{eq:loss_cfm}) were again chosen as simple feed-forward neural networks with Softplus activation functions. 
Networks used consisted of six layers with 500 neurons each. Models were trained for 1000 epochs with a batch size of 50 and a learning rate of 0.001. 
\ \\
All models were trained on up to 100\,000 samples of the respective target distributions. 
The trained models were then used to map the sparse grid quadrature nodes in both 9 and 25 dimensions to the target distributions. 
As an example, the mapped sparse grid nodes obtained from the ACF model for the bigamma distribution using $M=9$ modes are shown in \cref{fig:Sparse_Grid_INN}.

\subsection{Monomial tests}

As a first proof of concept, we tested the workflow outlined above on the integration of $M$-variate monomials. 
For each $k\in \mathbb N$ let $\mathcal{M}_k$ denote the set of all monomials in $M=9$ and $25$ variables of total degree less than or equal to $k$. 
As indicated in \cref{eq:rand_pol_error}, the sparse grid integration of each monomial $p\in\mathcal{M}_k$ is then carried out as the sum of the evaluations on the $n_\text{quad}$ transformed sparse grid quadrature nodes $\veta_j$ weighted by the sparse grid weights $w_j$. 
The value obtained is then compared to the mean of $p$ evaluated in $N=10^5$ MC samples $\veta^\text{MC}_i$ from the distribution defined in \cref{eq:LevyXi}, see \cref{app:DiscretizationLeyField} for the numerical implementation. 
Finally, the average error $\epsilon_k$ is calculated over all $p\in\mathcal{M}_k$.
\begin{align}\label{eq:rand_pol_error}
    \epsilon_k
    =
    \frac{1}{\lvert \mathcal{M}_k \rvert}
    \sum_{p\in \mathcal{M}_k}
    \left|
    \sum_{j=1}^{n_\text{quad}}
    w_j p(\veta_j) - \frac{1}{N}\sum_{i=1}^N p(\veta^\text{MC}_i)
    \right|.
\end{align}
These experiments were carried out for $M=9$ and $M=25$ variables for total degrees $k=1,2,3$, for which we employed the quadrature levels $L_\text{SG}=2,3$ and $4$, respectively.
As described in \cref{sec:SG}, these sparse grid quadrature rules are exact up to total degrees $3$, $5$ and $7$, respectively, so that, for an exact transport map, the quadrature error should be zero.
The number of quadrature nodes and the degree of exactness for the GH Smolyak quadrature rules used are summarized in \cref{tab:monimials}.

\begin{table}[ht]
\begin{center}
\begin{tabular}{|c|r|r|c|} \hline
$L_\text{SG}$ 
& \vtop{\hbox{\strut \# nodes }\hbox{\strut ($M=9$) }}
& \vtop{\hbox{\strut \# nodes }\hbox{\strut ($M=25$) }}
& \vtop{\hbox{\strut degree of }\hbox{\strut exactness }}\\ \hline
1  &     1 &      1 & 1\\
2  &    19 &     51 & 3\\
3  &   181 &  1 301 & 5 \\
4  & 1 177 & 22 201 & 7\\ \hline
\end{tabular}
\vspace{1em}
\end{center}
\caption{Summary of employed Gauss--Hermite Smolyak sparse quadrature levels $L_\text{SG}$, the associated number of quadrature nodes and degree of multivariate polynomial exactness.
\label{tab:monimials}
}
\end{table}

\subsection{Random diffusion problem}

As a specific instance of the stationary diffusion problem \cref{eq:diff_eq} we consider what is sometimes called the \emph{flow cell}, which models Darcy flow through an inhomogeneous medium on the unit square $D=[-1,1]^2$ with no sinks or sources, i.e.\ $f_\text{source}\equiv 0$, with homogeneous Neumann conditions at the upper and lower boundaries $\partial_{N} = [-1,1]\times \{-1, 1\}$ as well as a pressure difference driving the flow imposed by Dirichlet boundary conditions $g_D(x,y) = \frac{1}{2}(1-x)$ at the left and right $\partial_D = \{-1,1\}\times[-1,1]$. 

The equation is solved using a mixed finite element discretisation based on a triangular mesh with a piecewise constant representation of $u$ combined with a lowest-order Raviart--Thomas discretisation of the flux. 
We use three different meshes of 48, 798 and 8510 triangles, hereafter referred to as coarse, medium, and fine mesh. 
The goal of the computation is to approximate a scalar QoI given by the expectation of the flux through the right boundary.
\begin{figure}[!ht]
\includegraphics[width=.32\textwidth]{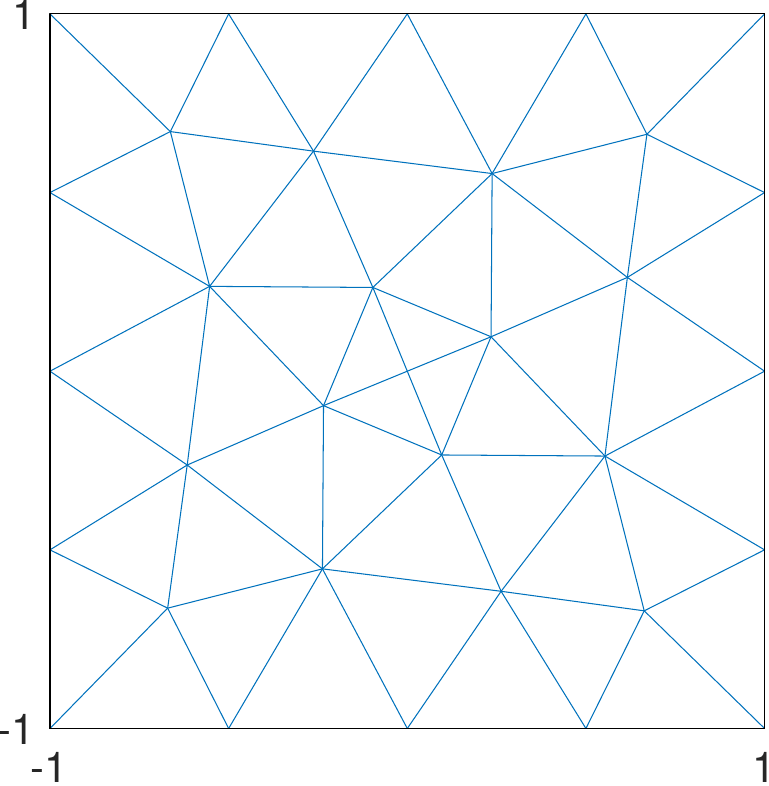}
\includegraphics[width=.32\textwidth]{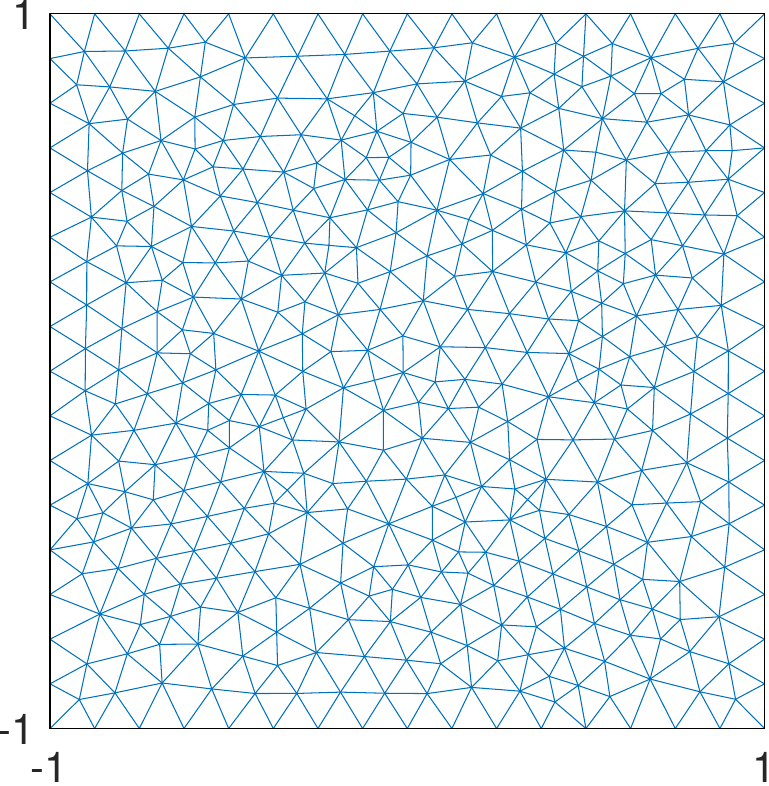}
\includegraphics[width=.32\textwidth]{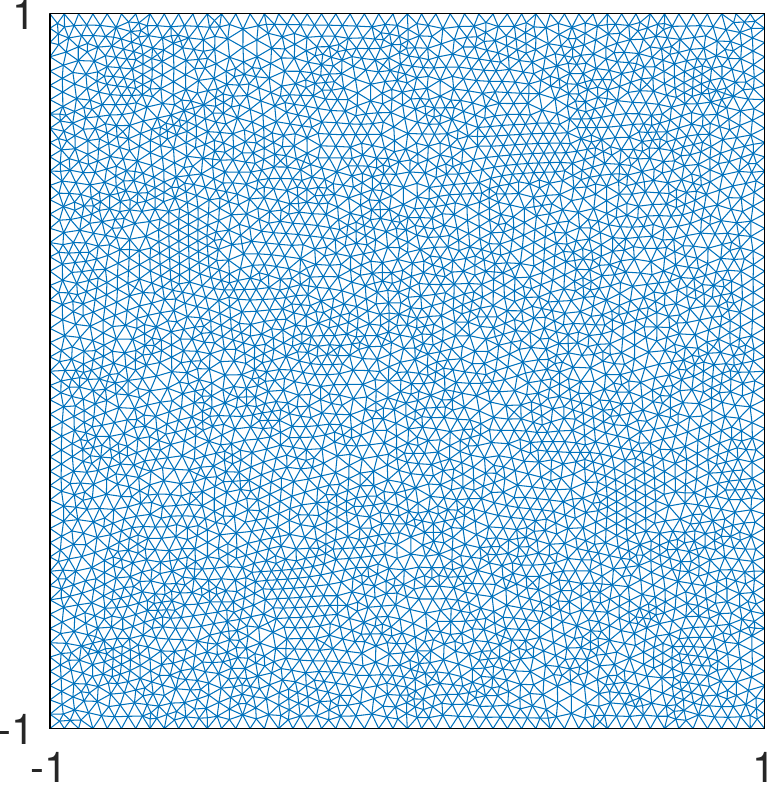}
\caption{The three meshes employed in the finite element solution consisting of 48, 798, and 8510 elements and referred to in the following as \emph{coarse}, \emph{medium} and \emph{fine}, respectively.}
\label{fig:meshes}
\end{figure}

The random conductivity $a$ is represented by a smoothed Lévy random field $Z_k$, for which we consider four common special cases, namely Gaussian, Poisson, gamma and bigamma noise fields convolved with a stationary Matérn covariance kernel with correlation length $m=0.1$ and smoothness parameter $\alpha=3$.
For each distribution we set the mean to zero and the variance to $\sigma^2 = 0.5$.

Since Gaussian and bigamma distributions take negative values, we applied the exponential function to all random fields. This results in positive realizations.
Poisson and gamma fields are always positive, so an exponential transformation is not needed for positivity. 
It is important to note that Poisson fields are always zero with positive probability, so we need to add a small constant to all Poisson fields to get samples that are bounded away from zero. 
However, we have used the exponential transformation in all cases for better comparability. 

\subsection{Discussion of error sources}

In the following section, we will discuss various error sources for both random monomials and the QoI of PDE solutions. 
We begin with two deterministic discretization aspects. 
We evaluate fixed realizations and compare different FE and Mercer approximations. 
Afterwards, we look at the accuracy of learned quadrature rules, where we used a MC estimator with $10^5$ samples as a reference solution. 
In our plots we indicate the level of the $95\%$ confidence interval of the corresponding MC estimate. 
We consider a quadrature successful if the error lies within the confidence interval. 
We will discuss the effect of the quadrature level and the regularity of the integrand.
Next, we compare results for four different distributions as well as three models: affine coupling flows (ACF), conditional flow matching (CFM) and optimal transport conditional flow matching (OTCFM). 
Finally, we discuss the quality of models for different amounts of training data.

\subsubsection{Finite element discretization error}

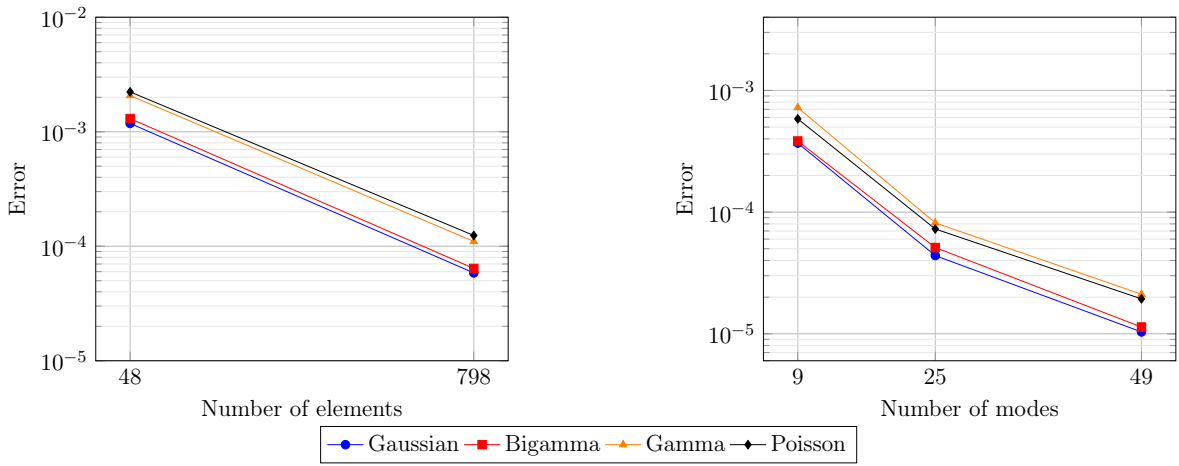
\begin{figure}[!ht]
  \centering
  \begin{tikzpicture}[scale=0.8]
      \begin{axis}[xmode=log, ymode=log,xlabel=Number of elements ,ylabel=Error,ymin=1e-5, ymax= 1e-2,legend to name = legenddiscerror, legend columns=-1, xtick={48,798},xticklabels={48,798}, grid=both,major grid style={black!25},minor grid style={black!10}]
	\addplot[blue, mark=*] table [x=elements, y=norm, col sep=comma]{images/FEMError/FEMesh.csv};
	\addlegendentry{Gaussian}
	\addplot[red, mark=square*] table [x=elements, y=bigamma, col sep=comma]{images/FEMError/FEMesh.csv};
	\addlegendentry{Bigamma}
	\addplot[orange, mark=triangle*] table [x=elements, y=gamma, col sep=comma]{images/FEMError/FEMesh.csv};
	\addlegendentry{Gamma}
	\addplot[black, mark=diamond*] table [x=elements, y=pois, col sep=comma]{images/FEMError/FEMesh.csv};
	\addlegendentry{Poisson}
      \end{axis}
    \end{tikzpicture}
  \hfill
    \begin{tikzpicture}[scale=0.8]
      \begin{axis}[ymode=log,xmode=linear,xlabel=Number of modes,ylabel=Error,ymin=6e-6, ymax=4e-3, xtick={9,25,49},xticklabels={9,25,49}, grid=both,major grid style={black!25},minor grid style={black!10}]
	\addplot[blue, mark=*] table [x=degree, y=norm, col sep=comma]{images/FEMError/FEDeg.csv};
	\addplot[red, mark=square*] table [x=degree, y=bigamma, col sep=comma]{images/FEMError/FEDeg.csv};
	\addplot[orange,mark=triangle*] table [x=degree, y=gamma, col sep=comma]{images/FEMError/FEDeg.csv};
	\addplot[black, mark=diamond*] table [x=degree, y=pois, col sep=comma]{images/FEMError/FEDeg.csv};
      \end{axis}
    \end{tikzpicture}

    \begin{center}
  \scalebox{0.8}{
    \pgfplotslegendfromname{legenddiscerror}
  }
\end{center}
	\caption{Discretization and truncation errors of the QoI for a random realization of the conductivity following four different probability distributions. 
    Left: error against finite element mesh size (full field without truncation). 
    Right: error against modes retained in Mercer series (fine mesh).}
  \label{fig:DiscError}
\end{figure}

Since we use the finite element method to solve the PDE, we do not know the exact solution, but only a finite dimensional approximation. 
As mentioned above, we work with three different meshes consisting of 48, 798 and 8510 elements, see \cref{fig:meshes}. 
The right-hand panel in \cref{fig:DiscError} shows how the discretisation error decays on a finer mesh. Here we have compared the solution obtained on the coarse and intermediate meshes with that on the fine mesh. For each realisation we used a sample of an untruncated random field.
The FE solution converges at a rate of 1 with respect to the number of elements in the mesh, which translates into a rate of $\mathcal{O}(h^2)$, where $h$ is the grid parameter, which is the expected result. Furthermore, the error is quite deterministic and behaves the same for all distributions considered.

\subsubsection{Mercer series truncation error}

The trained models do not operate on the full random fields but on a finite dimensional approximation. 
Specifically, we truncate the Mercer series after $(2n+1)^2$ terms (the constant and $n$ pairs of sine and cosine eigenfunctions in each direction, cf.\ \cref{app:DiscretizationLeyField}). 
The left panel of \cref{fig:DiscError}, where we compare the value of the quantity of interest of a random solution for each distribution, shows how the truncation error decreases as more terms are retained in the series.
As a reference, we used the (almost) untruncated random field obtained by discretising the torus into $100 \times 100$ squares, resulting in $100^2$ modes.
We observe that the truncation error is roughly the same for all distributions.

\subsubsection{Quadrature error}

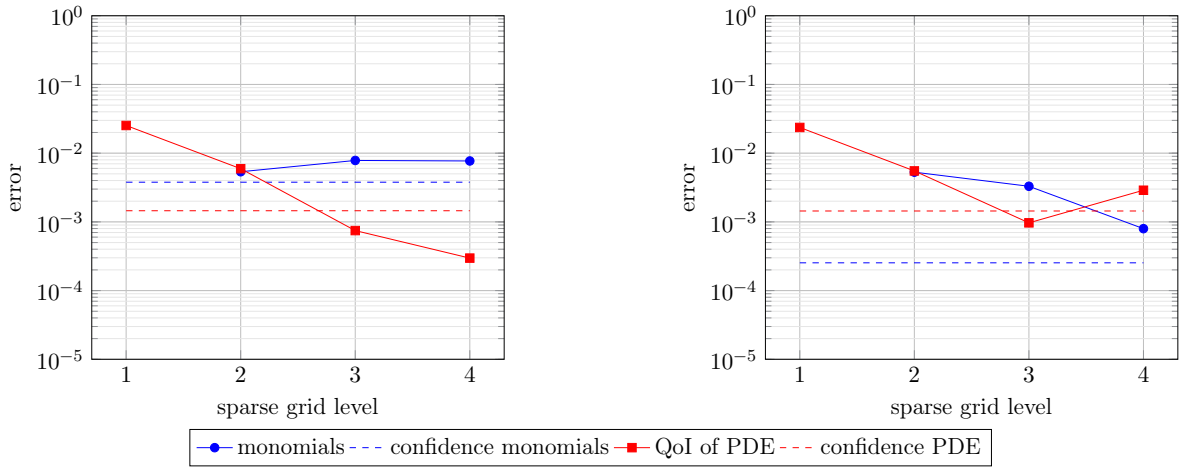
\begin{figure}[!ht]
    \begin{tikzpicture}[scale=0.8]
	  \begin{axis}[ymode=log,xlabel=sparse grid level ,ylabel=error,ymin=1e-5, ymax=1e0,legend to name=legendquaderror, legend columns=-1, xtick={1,2,3,4}, grid=both,major grid style={black!25},minor grid style={black!10}]
        
        \addplot[blue, mark=*] table [x=level, y=100000, col sep=comma]{images/Random_Polynomials/Pol_Errors_LevelsA2_OTCFM_bigamma.csv};
        \addlegendentry{monomials}
        \addplot[blue,dashed] table [x=level, y=bigamma, col sep=comma]{images/Random_Polynomials/CI_95_monomialsA2.csv};
        \addlegendentry{confidence monomials}
		\addplot[red, mark=square*] table [x=level, y=100000, col sep=comma]{images/FM_Models_03032025/ErrorDegCFMbigammaF.csv};
        \addlegendentry{QoI of PDE}
        \addplot[red,dashed] table [x=level, y=conf, col sep=comma]{images/FM_Models_03032025/ErrorDegCFMbigammaF.csv};
        \addlegendentry{confidence PDE}
	 \end{axis}
    \end{tikzpicture}
  \hfill
    \begin{tikzpicture}[scale=0.8]
	  \begin{axis}[ymode=log,xlabel=sparse grid level ,ylabel=error,ymin=1e-5, ymax=1e0, xtick={1,2,3,4}, grid=both,major grid style={black!25},minor grid style={black!10}]
        
        \addplot[blue, mark=*] table [x=level, y=100000, col sep=comma]{images/Random_Polynomials/Pol_Errors_LevelsA3_OTCFM_bigamma.csv};
        \addplot[blue,dashed] table [x=level, y=bigamma, col sep=comma]{images/Random_Polynomials/CI_95_monomialsA3.csv};
		\addplot[red, mark=square*] table [x=level, y=100000, col sep=comma]{images/FM_Models_03032025/ErrorDegOTCFMbigammaFA3.csv};
        \addplot[red,dashed] table [x=level, y=conf, col sep=comma]{images/FM_Models_03032025/ErrorDegOTCFMbigammaFA3.csv};
	 \end{axis}
    \end{tikzpicture}

    \begin{center}
  \scalebox{0.8}{
    \pgfplotslegendfromname{legendquaderror}
  }
\end{center}
  \caption{Quadrature error comparison in the bigamma case for monomials as well as for the quantity of interest of random PDEs for the best models, for varying degrees and largest training data size. Left: 9 modes. Right: 25 modes.}  
  \label{fig:rand_pol_results_quadrature}
\end{figure}

When integrating a function using quadrature, only polynomials up to a certain degree can be integrated exactly. 
This degree of exactness increases with the number of quadrature points used. \Cref{fig:rand_pol_results_quadrature} shows the decreasing quadrature error for increasing sparse quadrature levels for the best models with the most training data for the bigamma case. 
We used the OTCFM model with both 9 and 25 modes for the monomials and for the QoI of the PDE the 9-mode model was generated using CFM, while for the 25-mode model we used OTCFM, as these displayed the best performance. 
We compare the result of the quadrature for both experiments to a MC estimator using $10^5$ samples. 
In our plots we indicate the level of the corresponding $95\%$ confidence interval. 
We observe that the quadrature error for the PDE solution is mostly decreasing monotonically as the sparse grid level increases. 
However, for the monomials, the quadrature error does not decrease for the 9-mode model. Despite the increasing accuracy of the SG quadrature rule, the problem we are trying to solve becomes harder: With increasing level, monomials of higher degrees can be integrated exactly. 
After application of the transport map, these do not remain monomials of the same degree, most likely not even polynomials. 
Thus, the quadrature rule can no longer integrate the function exactly.

\subsubsection{Regularity constraints for sparse quadrature}

While sparse grid quadrature rules can converge extremely fast when the function is smooth, convergence can become slow or even impossible when regularity is lacking. As discussed earlier, the more the target distribution differs from the normal distribution, the more complex the transport map becomes and the more regularity is lost making it harder to find a suitable quadrature rule. In addition, the structure of the neural networks is important, especially the activation function. We have used the softplus activation function because it is smooth without kinks, as it maintains smoothness in each layer. Attempts to learn using activation functions with kinks, such as ReLU functions resulted in rather inaccurate quadratures. This can be taken as an indication that the requirements of the SG quadrature rules on the regularity of the integrand should be respected by the neural network architecture, see the remarks at the end of \cref{sec:SG} and equation \cref{eq:SG_Error}. 

\subsubsection{Mathematical assumptions on the distributions}

Next, we look at the distribution type as a source of error. Recall that all our trained models attempt to learn a continuous transport map that transforms the density of a standard normal distribution into the density of our target distribution. For random fields following a normal distribution \cref{eq:LevyLaws_gauss}, such a map is fairly easy to find. All we really need to do is scale each dimension by its corresponding variance. Even the bigamma distribution \cref{eq:LevyLaws_bigamma} behaves nicely: It takes values over the entire real line just like the normal distribution. This also applies to the random vector containing the random variables in the Mercer expansion. The modes for both distributions can take any real values. It is quite possible to find a smooth transport map without too much curvature.

Things become more complicated for the gamma distribution \cref{eq:LevyLaws_gamma}, which is still continuous, but only takes positive values. Unfortunately, this also results in a different behavior of the random mode coefficient vector. Its first component, i.e., the coefficient of the first eigenfunction, which is constant, must be positive to ensure the positivity of the random field. The subsequent components, belonging to sine and cosine functions, can have any sign, but are somewhat bounded by the previous components, as their amplitudes must be small enough to keep the random field positive. So we are looking for a continuous map that transforms values in $\R^n$ into a convex cone with only positive values in the first dimension. Such a transformation is impossible with a bijective map, resulting in extremely nonlinear transformations, which usually makes sparse grid quadratures less accurate.

Distributions \cref{eq:LevyLaws_pois} like the Poisson distribution are even harder to handle. Not only do they take only positive values, but they also have jumps, making the distribution partially discrete. Since bijective maps like ACF, CFM or OTCFM can only map a continuous source distribution to a continuous target distribution, the approximation of the discrete characteristics again requires a lot of capacity and high nonlinearity of the networks, resulting in poorly performing transport maps and thus rather inaccurate quadrature rules.

We display the quadrature errors of all trained models for each distribution in \cref{fig:ModelCompDeg03032025FA2}. We used the highest quadrature level on the finest grid for 9 modes. For the normal and Poisson distributions, the model trained via ACF performed best. For gamma and bigamma, we found the CFM model to be the best. As a reference solution, we used a MC estimator with $10^5$ samples. In the plots we indicate the level of its $95\%$ confidence interval.

\subsubsection{Comparison of different models}

\begin{figure}[!ht]
    \subfloat[Normal distribution]{\begin{tikzpicture}[scale=0.8]
	  \begin{axis}[xtick={1,2,3,4},ymode=log,xlabel=sparse grid level,ylabel=error,ymin=1e-5, ymax=1e0, legend columns=-1,  legend to name = legenddistributions,  grid=both,major grid style={black!25},minor grid style={black!10}]
        \addplot[blue, mark=*] table [x=level, y=100000, col sep=comma]{images/FM_Models_03032025/ErrorDegINNnormF.csv};
        \addlegendentry{ACF}
        \addplot[red, mark=square*] table [x=level, y=100000, col sep=comma]{images/FM_Models_03032025/ErrorDegCFMnormF.csv};
        \addlegendentry{CFM}
        \addplot[orange, mark=triangle*] table [x=level, y=100000, col sep=comma]{images/FM_Models_03032025/ErrorDegOTCFMnormF.csv};
        \addlegendentry{OTCFM}
    \addplot[black, dashed] table [x=level, y=conf, col sep=comma]{images/FM_Models_03032025/ErrorDegINNnormF.csv};
        \addlegendentry{confidence}
        \end{axis}
    \end{tikzpicture}}
  \hfill
    \subfloat[Bigamma distribution]{\begin{tikzpicture}[scale=0.8]
	  \begin{axis}[xtick={1,2,3,4}, ymode=log,xlabel=sparse grid level,ylabel=error,ymin=1e-5, ymax=1e0, grid=both,major grid style={black!25},minor grid style={black!10}]
        \addplot[blue, mark=*] table [x=level, y=100000, col sep=comma]{images/FM_Models_03032025/ErrorDegINNbigammaF.csv};
        \addplot[red, mark=square*] table [x=level, y=100000, col sep=comma]{images/FM_Models_03032025/ErrorDegCFMbigammaF.csv};
        \addplot[orange, mark=triangle*] table [x=level, y=100000, col sep=comma]{images/FM_Models_03032025/ErrorDegOTCFMbigammaF.csv};
    \addplot[black,dashed] table [x=level, y=conf, col sep=comma]{images/FM_Models_03032025/ErrorDegCFMbigammaF.csv};
        \end{axis}
    \end{tikzpicture}}
    \newline
    \subfloat[Gamma distribution]{\begin{tikzpicture}[scale=0.8]
	  \begin{axis}[xtick={1,2,3,4}, ymode=log,xlabel=sparse grid level,ylabel=error,ymin=1e-5, ymax=1e0, grid=both,major grid style={black!25},minor grid style={black!10}]
        \addplot[blue, mark=*] table [x=level, y=100000, col sep=comma]{images/FM_Models_03032025/ErrorDegINNgammaF.csv};
        \addplot[red, mark=square*] table [x=level, y=100000, col sep=comma]{images/FM_Models_03032025/ErrorDegCFMgammaF.csv};
        \addplot[orange, mark=triangle*] table [x=level, y=100000, col sep=comma]{images/FM_Models_03032025/ErrorDegOTCFMgammaF.csv};
    \addplot[black, dashed] table [x=level, y=conf, col sep=comma]{images/FM_Models_03032025/ErrorDegCFMgammaF.csv};
        \end{axis}
    \end{tikzpicture}}
  \hfill 
    \subfloat[Poisson distribution]{\begin{tikzpicture}[scale=0.8]
	  \begin{axis}[xtick={1,2,3,4}, ymode=log,xlabel=sparse grid level,ylabel=error,ymin=1e-5, ymax=1e0, grid=both,major grid style={black!25},minor grid style={black!10}]
        \addplot[blue, mark=*] table [x=level, y=100000, col sep=comma]{images/FM_Models_03032025/ErrorDegINNpoisF.csv};
		\addplot[red, mark=square*] table [x=level, y=100000, col sep=comma]{images/FM_Models_03032025/ErrorDegCFMpoisF.csv};
        \addplot[orange,mark=triangle*] table [x=level, y=100000, col sep=comma]{images/FM_Models_03032025/ErrorDegOTCFMpoisF.csv};
    \addplot[black, dashed] table [x=level, y=conf, col sep=comma]{images/FM_Models_03032025/ErrorDegINNpoisF.csv};
        \end{axis}
    \end{tikzpicture}}
\begin{center}
  \scalebox{0.8}{
    \pgfplotslegendfromname{legenddistributions}
  }
\end{center}

  \caption{Quadrature error of sparse quadrature rule learned by three different methods (9 modes). We compared the model with the highest amount of training data to a MC estimator with $10^5$ samples on the finest mesh.}   
 \label{fig:ModelCompDeg03032025FA2}
\end{figure}
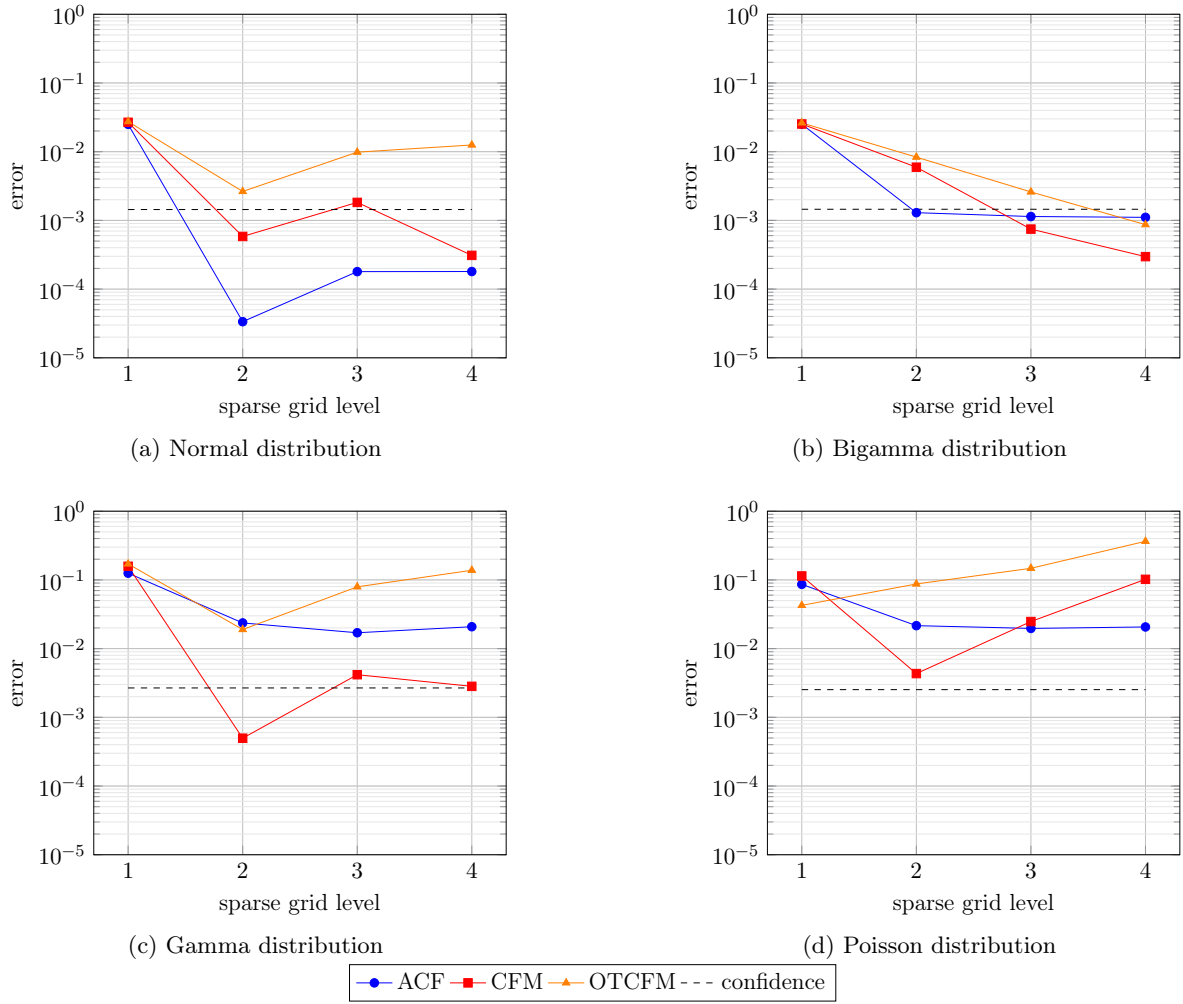

In these experiments, we considered three different models for the learning process of the transport map: ACF, CFM and OTCFM. Their performance is shown in \cref{fig:ModelCompDeg03032025FA2}. For the normal distribution, the OTCFM method failed to learn a suitable quadrature rule, while the ACF gives the best results, even at lower quadrature levels. For the bigamma distribution, all three methods learned a quadrature rule with a result within the $95\%$ confidence interval of the MC estimator. For the ACF, the difference between levels 2 to 4 is rather small, while the other two methods slowly gain accuracy with increasing quadrature level. For gamma fields, only the quadrature rule learned by the CFM gave a result close to the confidence interval. Finally, for the Poisson case, all three methods failed as expected.

\subsubsection{Variance of the training data and the optimizer}

\begin{figure}[!ht]
    \begin{tikzpicture}[scale=0.8] 
    \begin{axis}[xmode=log, ymode=log,xlabel=sample size,ylabel=error,ymin=1e-5, ymax=1e0,legend to name=legend_random_seeds, legend columns=-1, grid=both,major grid style={black!25},minor grid style={black!10},tick label style={font=\footnotesize},label style={font=\footnotesize},legend style={font=\footnotesize}]
	\addplot[name path=upper5, draw=none, forget plot] table[x=trainsize, y=4confU, col sep=comma] {images/var_training/ErrorTrainsizeCFMbigammaFA2.csv};
	\addplot[name path=lower5, draw=none, forget plot] table[x=trainsize, y=4confL, col sep=comma] {images/var_training/ErrorTrainsizeCFMbigammaFA2.csv};
	\addplot[fill=blue!20, fill opacity=.5, draw=none] fill between[of=upper5 and lower5];
	\addlegendentry{$2\sigma$ interval of error}
	\addplot[blue!60, thick, forget plot] table[x=trainsize, y=4confU, col sep=comma] {images/var_training/ErrorTrainsizeCFMbigammaFA2.csv};
	\addplot[blue!60, thick, forget plot] table[x=trainsize, y=4confL, col sep=comma] {images/var_training/ErrorTrainsizeCFMbigammaFA2.csv};
        \addplot[black, mark=*] table [x=trainsize, y=4, col sep=comma]{images/var_training/ErrorTrainsizeCFMbigammaFA2.csv};
	\addlegendentry{error mean}
        \addplot[red, dashed] table [x=trainsize, y=conf, col sep=comma]{images/var_training/ErrorTrainsizeCFMbigammaFA2.csv};
	\addlegendentry{tolerance of reference solution}
        \end{axis}
    \end{tikzpicture}
  \hfill 
    \begin{tikzpicture}[scale=0.8]
	  \begin{axis}[xmode=log,ymode=log,xlabel=sample size,ylabel=error,ymin=1e-5, ymax=1e0, grid=both,major grid style={black!25},minor grid style={black!10},tick label style={font=\footnotesize},label style={font=\footnotesize},legend style={font=\footnotesize}]
	\addplot[name path=upper6, draw=none] table[x=trainsize, y=4confU, col sep=comma] {images/var_training/ErrorTrainsizeOTCFMbigammaMA2.csv};
	\addplot[name path=lower6, draw=none] table[x=trainsize, y=4confL, col sep=comma] {images/var_training/ErrorTrainsizeOTCFMbigammaMA2.csv};
	\addplot[fill=blue!20, fill opacity=.5, draw=none] fill between[of=upper6 and lower6];
	\addplot[blue!60, thick] table[x=trainsize, y=4confU, col sep=comma] {images/var_training/ErrorTrainsizeOTCFMbigammaMA2.csv};
	\addplot[blue!60, thick] table[x=trainsize, y=4confL, col sep=comma] {images/var_training/ErrorTrainsizeOTCFMbigammaMA2.csv};
        \addplot[black, mark=*] table [x=trainsize, y=4, col sep=comma]{images/var_training/ErrorTrainsizeOTCFMbigammaMA2.csv};
        \addplot[red, dashed] table [x=trainsize, y=conf, col sep=comma]{images/var_training/ErrorTrainsizeOTCFMbigammaMA2.csv};
        \end{axis}
    \end{tikzpicture}
\begin{center}
  \scalebox{0.8}{
    \pgfplotslegendfromname{legend_random_seeds}
  }
\end{center}
  \caption{Mean and standard deviation of learned quadrature rules for the 9-mode bigamma model using $10$ random seeds. Left: CFM. Right: OTCFM.}  
  \label{fig:random_seeds}
\end{figure}

\oge{
The training of neural networks highly depends on the training data used as well as the training process. 
The optimizer of our choice was the AdamW algorithm, which generates stochastic gradients from a random subset of the training data. 
In this way, the entire training process depends on a random seed. 
In order to measure the standard deviation of the quadrature error resulting from our method, we repeated the training procedure for $10$  different random seeds for the 9-mode bigamma distribution and evaluated the resulting quadrature rules. 
\Cref{fig:random_seeds} shows the mean of the quadrature error as well as the surrounding $2\sigma$-intervals as a measure for the random variability of the training process.
}

\subsubsection{Learning error}

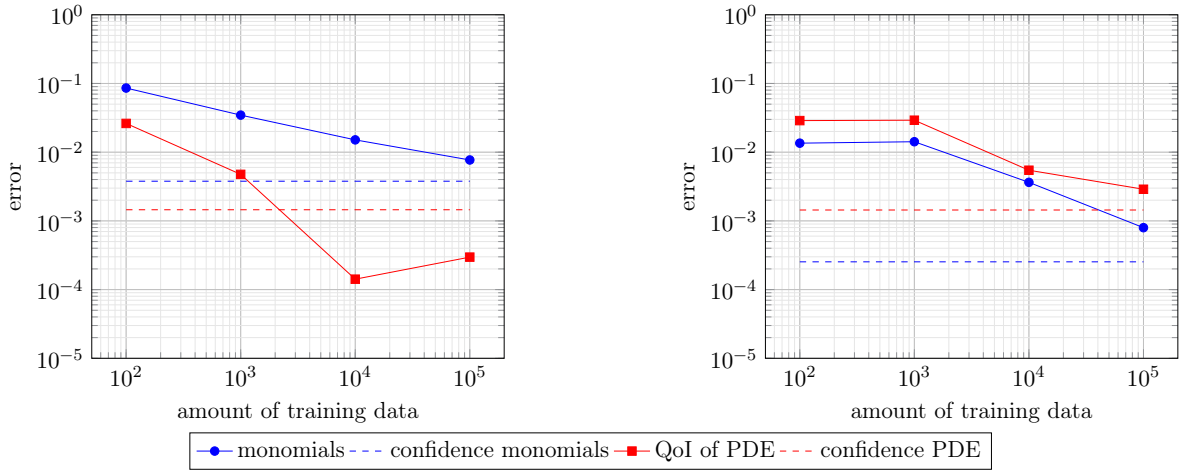
\begin{figure}[!ht]
    \begin{tikzpicture}[scale=0.8]
      \begin{axis}[xmode=log, ymode=log,xlabel=amount of training data,ylabel=error,ymin=1e-5, ymax=1e0,legend to name=legendlearningerror, legend columns=-1, grid=both,major grid style={black!25},minor grid style={black!10}]
        \addplot[blue, mark=*] table [x=Datasize, y=deg 3, col sep=comma]{images/Random_Polynomials/Pol_Errors_A2_OTCFM_bigamma.csv};
        \addlegendentry{monomials}
        \addplot[blue, dashed] table [x=trainsize, y=bigamma, col sep=comma]{images/Random_Polynomials/CI_95_monomialsA2.csv};
        \addlegendentry{confidence monomials}
        \addplot[red, mark=square*] table [x=trainsize, y=4, col sep=comma]{images/FM_Models_03032025/ErrorTrainsizeCFMbigammaF.csv};
        \addlegendentry{QoI of PDE}
        \addplot[red, dashed] table [x=trainsize, y=conf, col sep=comma]{images/FM_Models_03032025/ErrorTrainsizeCFMbigammaF.csv};
        \addlegendentry{confidence PDE}
        \end{axis}
    \end{tikzpicture}
  \hfill
    \begin{tikzpicture}[scale=0.8]
      \begin{axis}[xmode=log, ymode=log,xlabel=amount of training data,ylabel=error,ymin=1e-5, ymax=1e0, grid=both,major grid style={black!25},minor grid style={black!10}]
        \addplot[blue, mark=*] table [x=Datasize, y=deg 3, col sep=comma]{images/Random_Polynomials/Pol_Errors_A3_OTCFM_bigamma.csv};
        \addplot[blue, dashed] table [x=trainsize, y=bigamma, col sep=comma]{images/Random_Polynomials/CI_95_monomialsA3.csv};
        \addplot[red, mark=square*] table [x=trainsize, y=4, col sep=comma]{images/FM_Models_03032025/ErrorTrainsizeOTCFMbigammaFA3.csv};
        \addplot[red, dashed] table [x=trainsize, y=conf, col sep=comma]{images/FM_Models_03032025/ErrorTrainsizeOTCFMbigammaFA3.csv};
        \end{axis}
    \end{tikzpicture}
    \begin{center}
  \scalebox{0.8}{
    \pgfplotslegendfromname{legendlearningerror}
  }
\end{center}
  \caption{Statistical error comparison in the bigamma case for monomials as well as for the quantity of interest of random PDEs for the best models for varying training sizes and highest quadrature degree. Left: 9 modes. Right: 25 modes.}  
  \label{fig:rand_pol_results_learning}
\end{figure}

Finally, we consider the statistical learning error that results from the approximation of the true transport map by a machine learning model. The goal of all the generative models discussed so far is to find a suitable transport map that transforms a multivariate standard normal distribution into the desired distribution of the random vector containing the random variables in the Mercer expansion. We want the learned distribution to be as close as possible to the true distribution. By increasing the amount of training data, we can train longer without overfitting, resulting in more accurate models.
\Cref{fig:rand_pol_results_learning} shows the error reduction for increasing training data sizes for the integration of monomials and the QoI from our simulations, respectively. Here we used the best model for the bigamma case on sparse grid level 4, which is the highest level we used, which is OTCFM. For the integration of monomials, we used the model trained via OTCFM with both 9 and 25 modes. For the QoI of the diffusion equation we used the CFM solution for 9 modes and the OTCFM solution for 25 modes.

\subsubsection{Comparison with Monte Carlo baseline}
\begin{figure}[!ht]
    \begin{tikzpicture}[scale=0.8]
      \begin{axis}[xmode=log, ymode=log,xlabel=number of function evaluations,ylabel=error,ymin=1e-5, ymax=1e0,legend to name=legend_mmd, legend columns=-1, grid=both,major grid style={black!25},minor grid style={black!10}]
        \addplot[blue] table [x=samples, y=mean, col sep=comma]{images/baselines/ErrorMCA2.csv};
        \addlegendentry{Monte Carlo}
        \addplot[orange,mark=triangle*] table [x=samples, y=100000, col sep=comma]{images/var_training/ErrorDegCFMbigammaFA2.csv};
        \addlegendentry{CFM}
	\addplot[black, mark=diamond*] table [x=samples, y=100000, col sep=comma]{images/var_training/ErrorDegOTCFMbigammaMA2.csv};
    \addlegendentry{OTCFM}
        \end{axis}
    \end{tikzpicture}
  \hfill
    \begin{tikzpicture}[scale=0.8]
      \begin{axis}[xmode=log, ymode=log,xlabel=number of function evaluations,ylabel=standard deviation of error,ymin=1e-5, ymax=1e0, legend columns=-1, grid=both,major grid style={black!25},minor grid style={black!10}]
        \addplot[blue] table [x=samples, y=std, col sep=comma]{images/baselines/ErrorMCA2.csv};
        \addplot[orange,mark=triangle*] table [x=samples, y=std, col sep=comma]{images/var_training/ErrorDegCFMbigammaFA2.csv};
	\addplot[black, mark=diamond*] table [x=samples, y=std, col sep=comma]{images/var_training/ErrorDegOTCFMbigammaMA2.csv};
        \end{axis}
    \end{tikzpicture}
\begin{center}
  \scalebox{0.8}{
    \pgfplotslegendfromname{legend_mmd}
  }
\end{center}
  \caption{Comparison of sample efficiency for plain MC and quadrature learned by (OT-)CFM. 
  Left: mean error. 
  Right: standard deviation of error.}  
  \label{fig:baselines}
\end{figure}

\oge{
To compare our results with a plain MC baseline, we consider $10^5$ randomly generated samples from the 9-mode bigamma distribution. As an analytical evaluation of the QoI integrals with respect to the resulting L\'evy distribution is not possible, we take the average of the QoI over the $10^5$ samples as a reference value. 
We estimate the mean squared error of the MC estimates via non-parametric bootstrapping \cite{efron1994introduction}: we draw a number $n=1,\ldots,10^5$ of random samples with replacement from the  original $10^5$ samples and compute the mean QoI for the sampled data. 
We then compute the absolute difference between the original mean QoI for the full data set and the sample to compute the error. 
For each value of $n$, this procedure is repeated for $1\,000$ times and the error averaged. 
This estimates the mean error of the MC approach with $n$ evaluations to high accuracy. 
The blue solid line in \Cref{fig:baselines} (left) shows the resulting errors as a function of the number of required PDE solves $n$. 
We recognize the $n^{-1/2}$ scaling of the central limit theorem.
}

\oge{
For comparison, we train our CFM and OTCFM models on the full data set of size $10^5$ containing the modes. 
Note that no PDE solves are required in this step. This training is repeated $10$ times for $10$ different random initializations. We then evaluate each of the $10$ CFM or OTCFM models on sparse grids of increasing level $k=1,\ldots,4$ and compute the corresponding finite element approximations of the QoI and then plot the average absolute error for the QoI integrals relative to the reference MC solution as a function of the required PDE solves. 
The resulting graphs are the orange (CFM) and black (OTCFM) lines Figure \ref{fig:baselines} (left). 
We clearly see that the average MC error lies significantly above both the average CFM and OTCFM errors.}

\oge{
We also compare the robustness of MC for computing QoI expectations with the learned CFM and OTCFM sparse grid quadrature by computing standard deviations over  $1\,000$ repetitions of bootstrap sampling and the $10$ repetitions of neural network training, see \Cref{fig:baselines} (right). 
The standard deviations of CFM and OTCFM sparse grid quadratures are considerably lower than those of the MC method, as sparse grid quadrature nodes are non-random. 
In both application-relevant metrics, the average precision and the robustness, our algorithm thus outperforms  MC sampling with the same number of QoI evaluations. }

\section{Conclusions and outlook} \label{sec:conclusions}
 
We have discussed learned sparse grid quadratures for complex distributions for which the analytic form of the probability density function is unknown.
At the same time, sampling from the complex distribution is assumed to be inexpensive compared to the evaluation of the quantity of interest.

In this situation, we propose learning sparse grid quadrature rules by normalizing flows mapping the distribution at hand (approximately) to a standard normal distribution. Normalizing flow maps are learned from the available data and then inverted and evaluated on sparse grid quadrature points. In this way, we map Gaussian sparse grid quadratures to (approximate) sparse quadratures of the underlying distributions.

Three different normalizing flow architectures were tested -- affine coupling flows, flow matching and optimal transport flow matching. 
The results are compared to naive MC estimation with $10^5$ samples. 
We have considered two test cases: 
On the one hand, we assess the quadrature errors for multivariate monomials well within the exactness range of the unmapped Gauss--Hermite Smolyak rules. 
Here we observe a clear trend towards better numerical integration results when the sample size increases.

Our main example is the diffusion equation with a (smoothed and exponentiated) L\'evy random field as conductivity coefficient. 
Unlike in the Gaussian case, the L\'evy distribution cannot be mapped to a standard normal distribution by a simple linear transformation. 
Nevertheless, a modal expansion resembling the well known Karhunen-Loève expansion in the Gaussian case is possible \cite{ErnstEtAl2021TR,Reese2021}. 
We therefore take this as a suitable test case for our concept of learned nonlinear transformations by  normalizing flows and learned quadrature rules.       
Our empirical results indicate that the learned quadrature rules indeed approximate the expected values of the quantity of interest. 
For our analysis of the remaining error, we provide several convergence studies and compare three different  normalizing flows, three levels of finite element discretization of the stationary diffusion equation, two different truncation levels of the modal expansion, four different probability laws for Lévy random fields, up to six different levels of sparse grid quadrature, and varying amounts of training data. 

Although empirically our `Learning to Integrate' approach produced results mostly in or close to the 95\% confidence interval of a MC estimation with $10^5$ samples and outperformed the MC baseline, our results are still in an experimental stage.  
Testing on more diverse sets of distributions is required before a best practice can be established. 

The theoretical status of the `Learning to Integrate' approach needs to be clarified. 
As all normalizing flow models employed are recent, the statistical learning theory for these is not yet fully developed. 
This still hinders a rigorous numerical and statistical combined analysis that results in guarantees for convergence in the limit of large samples and adaptively enlarged neural network architectures. 
Initial results in this direction can be found in \cite{GottschalkEtAl2025-PP}.

\vspace{0.5cm}
\paragraph{Acknowledgments} 
The authors thank Thomas Kalmes, Emil Partow,  Marco Reese and Tobias Riedlinger  for interesting discussions, Klaus Ritter for pointing out recent references on Gauss--Hermite quadrature and Yoshihito Kazashi for sharing his results prior to publication.

\bibliographystyle{siamplain}
\bibliography{lti}

\newpage
\begin{appendices}

\crefalias{section}{appendix}
\crefalias{subsection}{appendix}
\crefname{appendix}{Appendix}{Appendices}
\Crefname{appendix}{Appendix}{Appendices}
\section{Supplementary Material} 
In this appendix, we provide details on the numerical discretization and simulation of the L\'evy random fields and the solution of the diffusion equation for a given random field as coefficient.

\subsection{Discretization and fast simulation of random fields}
\label{app:DiscretizationLeyField}
In this appendix we present the discretization and generalized Mercer expansion \cite{ErnstEtAl2021TR} applied to the L\'evy random fields after circular embedding.

We have to discretize the torus for an efficient simulation of the Lévy random field coefficients with respect to eigenfunctions of $(-\Delta+m^2)^{-\alpha}$ that provides the training data for the normalizing flow. 
Let us, for the sake of concreteness, rescale the length scale of our computational domain so that $\mathbb{T}^d$ can be identified with $[-1,1]^d$ with identified opposing boundaries. 
The eigenfunctions of $(-\Delta+m^2)^{-\alpha}$ are the Fourier modes $\e^{i \kappa\cdot x}$ with $\kappa\in \pi\mathbb{Z}^d$ and $i$ the imaginary unit. 
The associated eigenvalues for the pseudodifferential operator $(-\Delta+m^2)^{-\alpha}$  are then $(|\kappa|^2+m^2)^{-\alpha}$.  
As the Mercer expansion of the Mat\'ern kernels requires real-valued eigenfunctions, we obtain these  by taking the real part $\cos( \kappa\cdot x)$ and the imaginary part $\sin(\kappa\cdot x)$ of $\e^{i \kappa\cdot x}$. 
It can be easily verified that these are real-valued eigenfunctions of $(-\Delta+m^2)^{-\alpha}$ for the same eigenvalue as $\e^{i \kappa\cdot x}$. 
Note that these functions are not normalized; the normalization constants are 
$2^{-\frac{d}{2}}$ for $\e^{i \kappa\cdot x}$ and $2^{\frac{1-d}{2}}$ for $\cos(\kappa\cdot x)$ and $\sin(\kappa\cdot x)$ for $\kappa\not=0$. 
The Mercer expansion of the smoothing operator in the given setting therefore coincides with the Fourier expansion of the (circulant) integral kernel $K^c_\alpha$ on $\mathbb{T}^d$ and is given by
\[
    K_{\alpha}^c(x-x')
    =
    \frac{1}{2^{d}}
    \sum_{\kappa\in\pi \mathbb{Z}^d}
    \frac{\e^{i\kappa\cdot x}\e^{-i\kappa\cdot x'}}
         {(|\kappa|^2+m^2)^{\alpha} }
    =
    \frac{1}{2^{d}}
    \sum_{\kappa\in \pi \mathbb{Z}^d}
    \frac{\cos(\kappa\cdot x)\cos(\kappa\cdot x')+\sin(\kappa\cdot x)\sin(\kappa\cdot x')}      {(|\kappa|^2+m^2)^{\alpha} }.
\]
The Mercer series is truncated by restricting the multi-index set $\mathbb{Z}^d$ to 
\[
    \mathbb{Z}^d_r = \{z\in\mathbb{Z}^d : |z|_\infty\leq r\}, \qquad r=0,1,2,\ldots,
\]
yielding an approximate expansion consisting of $M = |\mathbb{Z}^d_r|$ terms, which in turn results in an approximate random field given by
\begin{equation} \label{eqa:MercerMatern}
\begin{split}
    Z_{K_{\alpha,r}^c}(x)
    &=
    Z(K_{\alpha,r}^c(x-\cdot)) \\
    &=
    \frac{1}{2^{d}}
    \sum_{\kappa\in \pi \mathbb{Z}_r^d} 
    \frac{\cos(\kappa\cdot x)Z(\cos(\kappa\cdot (\cdot)))
          +\sin(\kappa\cdot x)Z(\sin(\kappa\cdot(\cdot) ))}
         {(|\kappa|^2+m^2)^{\alpha} }.
\end{split}         
\end{equation}

However, not all $2\times (2r+1)^d$ modes $Z(\sin(\kappa\cdot(\cdot)))$, $Z(\cos(\kappa\cdot(\cdot)))$, $\kappa\in\pi\mathbb{Z}^d_r$ are linearly independent. 
Indeed, for $\kappa=0$ we have $\sin(\kappa\cdot x' ) \equiv 0$ and consequently $Z(\sin(\kappa\cdot(\cdot) ))=0$. 
In addition, the symmetries $\cos(-\kappa \cdot x')=\cos(\kappa \cdot x')$ and $\sin(-\kappa \cdot x')=-\sin(\kappa \cdot x')$ and therefore $Z(\cos(-\kappa \cdot (\cdot)))=Z(\cos(\kappa \cdot(\cdot)))$ and $Z(\sin(-\kappa \cdot(\cdot)))=-Z(\sin(\kappa \cdot(\cdot)))$, further reduce the number of linearly independent modes to $M=(2r+1)^d$. 
Denoting by $\pi\mathbb{Z}_{r,+}^d$ those $\kappa=(\kappa_1,\ldots,\kappa_d) \in \pi\mathbb{Z}_{r}^d$ not identically zero for which the first nonzero component is positive, the approximation \cref{eqa:MercerMatern} yields the generalized Karhunen--Loève expansion for the random field $Z_{K_\alpha^c}$ as displayed in \eqref{eqa:MercerMatern-indep}.

In the next step we discretize the noise field, by averaging over a grid in space. To this aim, let $[-1,1]^q$ be subdivided in $\left(\frac{1}{L}\right)^q$ squares with edge length $2L$. Here we assume that $\frac{1}{L}=1,3,5\ldots$ is an odd integer. Let $\Gamma$ be the lattice of center points of these squares. The random value of the noise field with L\'evy characteristic $\psi(t)$ averaged over any such square $\Lambda$ are i.i.d.\ random variables with L\'evy distribution given by the characteristic function
\begin{equation} \label{eq:CharFunDisceteLevyNoise}
\ev{e^{\mathbbm{i}\, t Z(\mathbbm{1}_\Lambda)}}
=
\ev{e^{\mathbbm{i}\,  Z(t\mathbbm{1}_\Lambda)}}
= 
e^{\int_{\mathbb{T}^q}\psi\left(t\,\mathbbm{1}_\Lambda\right)\,\mathrm{d}x}
=
e^{(2L)^q\psi(t)},
\qquad t \in \mathbb R.
\end{equation}

Hence, a simulation of the noise field $Z$ on the lattice $\Gamma$ is obtained by drawing i.i.d. random variables with law determined by the above characteristic function. In this way, the discretized L\'evy field $Z_\Gamma$ is obtained, where $Z_\Gamma(x)$ gives the random value obtained for the square with center $x\in\Gamma$. 

On the given lattice, the lattice Laplacian $\Delta_\Gamma$ acts as
\[
\Delta_\Gamma h(x)=\left(\frac{1}{2L}\right)^2\sum_{y\in \Gamma_x} [h(y)-h(x)]=\left(\frac{1}{2L}\right)^2\left[\sum_{y\in \Gamma_x}h(y)-dh(x)\right],~~x\in \Gamma,
\]
where $\Gamma_x$ is the set of nearest neighbors of $x$ in $\Gamma$. On the space of functions $h:\Gamma\to \mathbb{C}$, the lattice Laplacian again is diagonalized by the 
eigen functions $e^{i\kappa\cdot x}$ where $\kappa\in \Gamma'=\pi\{-L,\ldots,L\}^q$ is the dual group to the discrete torus group $\Gamma$ \cite{van2012potential}, i.e. for $\kappa\in \Gamma'$, $e^{i\kappa\cdot (\cdot)}:\Gamma\to \{z\in\mathbb{C}: |z|=1\}$ is a character, i.e. a group homomorphism with respect to addition in $\Gamma$ and multiplication in $\mathbb{C}$. The eigen values  of these functions with respect to $\Delta_\Gamma$ are given by
\begin{equation}
\label{eq:SymbolDiscLaplacian}
\left(\frac{1}{2L}\right)^2\left(\sum_{j=1}^d \cos\left(2L \kappa_j\right)-d\right),~~\kappa\in \Gamma'.
\end{equation}

Thus, introducing the discrete Fourier transform $\mathcal{F}(h)(k)=\sum_{x\in\Gamma}h(x)e^{i\kappa\cdot x}$, $\kappa\in\Gamma'$, and its inverse $\mathcal{F}^{-1}$, we define 
\[
(-\Delta_\Gamma+m^2)^\alpha h(x)=\mathcal{F}^{-1}\left((-M_{\Gamma'}+m^2)^\alpha \mathcal{F}h\right)(x),~~x\in\Gamma,
\]
where $M_{\Gamma'}$ is the multiplication operator with the symbol of the discretized Laplacian given in \eqref{eq:SymbolDiscLaplacian}. It is easily seen that the operator $(-\Delta_\Gamma+m^2)^{-\alpha}$ is realized by   $(-\Delta_\Gamma+m^2)^{-\alpha} h(x)=\mathcal{F}^{-1}\left((-M_{\Gamma'}+m^2)^{-\alpha} \mathcal{F}(h)\right)(x),~~x\in\Gamma$. Consequently, the discretized random field $Z_{k,\Gamma}$ is obtained as the solution of the discretized stochastic pseudodifferential equation \eqref{eq:SPPDE},
\[
(-\Delta_\Gamma+m^2)^\alpha Z_{K,\Gamma}(x)=Z_\Gamma(x),~~x\in\Gamma ~~\Leftrightarrow Z_{K,\Gamma}(x)= [(-\Delta_\Gamma+m^2)^{-\alpha}Z_\Gamma](x),~~x\in\Gamma.
\]
This now gives an easy to implement and fast algorithm to simulate $Z_{K,\Gamma}$: First, for any of the $L^d$ squares $\Lambda(x)$ centred at $x\in\Gamma$, draw one i.i.d. value from the distribution with characteristic function \eqref{eq:CharFunDisceteLevyNoise} to obtain $Z_\Gamma(x)$. Second, apply the fast Fourier transform (FFT) implementation of the discrete Fourier transform to $Z_\Gamma$. Then, for any $k\in\Gamma'$, multiply the value obtained with $(-M_\Gamma+m^2)^{-\alpha}$. Fourth and lastly, apply the inverse FFT.    

The mode approximation of a random field now utilizes the fact that on the periodic lattice, the trigonometric functions $\cos(\pi \kappa\cdot x)$ and $\sin(\pi \kappa\cdot x)$ with $\kappa\in\Gamma'$ are the eigen-functions of $\Delta_\Gamma$. In order to provide a mode expansion, we only have to delete all those mode entries for the FFT of the noise, that should not be considered. Note that here we also utilize that the convolution operation with the kernel of $(-\Delta_\Gamma+m^2)^{-\alpha}$ acts as a multiplication operator in frequency space and therefore deleted modes remain deleted.

We pick the modes in the following way: As $L\to0$, the eigen values \eqref{eq:SymbolDiscLaplacian} converge to $(|k|^2+m^2)^{-\alpha}$. Thus, in the prescribed limit, they fall with the euclidean norm of the wave vector $k$. Therefore we pick those wave vectors $k$ closest to the origin. Here we work in $d=2$. Thus, in order to pick the largest eigenvalues in a symmetric manner, we select squared regions in $\Gamma'$ with edge length $1,3,5,\ldots$ containing $1,9,25\ldots$ modes with the highest eigen values. With the choice $1$, the random field just consists of the constant mode, which is of little interest. Therefore, we consider modes with $9$ and $25$ modes for the numerical experiments in the next section. 

The numerical implementation of the random fields has been performed with \texttt{R 4.5.0} using the built in \texttt{fft} command, which is a \texttt{C}-port of the original \texttt{Fortran} code in \cite{singleton1979mixed}. The sparse grid generation used the \texttt{SparseGrid} library version \texttt{0.8.2} by J.\ Ympa, which is an \texttt{R}-port of the \texttt{matlab} library from \cite{heiss2008likelihood}. For the interpolation to the FEA quadrature points, bi-linear interpolation with the \texttt{akima 0.6-3.4
}-package was utilized. The code will be made publicly available on \texttt{github}\footnote{\url{https://github.com/ToniKowalewitz/LearningToIntegrate}}.


\subsection{Mixed Finite Element Discretization} \label{app:MixedFE}

Since we are mainly interested in the resulting flow instead of the underlying potential, we employ the mixed formulation of the PDE, which introduces the flux $\vsigma = a\nabla u$ as an additional variable and solves the coupled first-order system
\begin{align}
  \vsigma - a\nabla u &= 0     \label{eq:mixed_coupling}\\
  -\nabla \cdot \vsigma &= f_\text{source}.  \label{eq:mixed_flux}
\end{align}

To obtain a variational formulation of our PDE, we multiply both \eqref{eq:mixed_coupling} and \eqref{eq:mixed_flux} by suitable test functions and integrate over the domain $D$.
Testing first \eqref{eq:mixed_flux} against $v\in L^2(D)$ gives
\begin{equation}
  -\int_D (\nabla\cdot\vsigma) v \,\d x = \int_D f_\text{source}v\,\d x\quad \text{for all }v\in L^2(D).
  \label{eq:var_flux}
\end{equation}
For \eqref{eq:mixed_coupling}, we use test functions 
$\vtau \in H_0(\text{div},D) := \{\vtau \in H(\text{div}; D) : \vtau \cdot \vn|_{\partial_N} = 0\}$ 
and divide by $a$, giving
\begin{equation*}
  \int_D \frac{1}{a}\vsigma \cdot \vtau\,\d x - \int_D \nabla u \cdot \vtau\,\d x = 0 
  \quad \text{for all }\vtau\in H_0(\text{div}; D).
\end{equation*}
We use the integration by parts formula
\begin{equation*}
  \int_D \nabla\cdot(u\vsigma) \,\d x 
  = 
  \int_D [\nabla u \cdot \vsigma + u (\nabla \cdot \vsigma)] \,\d x 
  = 
  \int_{\partial D} u \vsigma \cdot \vn \,\d s
\end{equation*}
and obtain
\begin{equation}  
  \int_D \frac{1}{a}\vsigma \cdot \vtau\,\d x + \int_D  u (\nabla \cdot \vtau) \,\d x 
  = 
  \int_{\partial_D} g_D \vtau\cdot \vn \,\d s \quad \text{for all }\vtau\in H_0(\text{div}; D).
  \label{eq:var_coupling}
\end{equation}
Note that only the Dirichlet part of the boundary conditions remains on the right hand side of the equation, as the homogeneous Neumann boundary condition is built into the solution space $H_0(\text{div}; D)$.
We can now combine \eqref{eq:var_coupling} and \eqref{eq:var_flux} again into
\begin{align*}
  A(\vsigma,\vtau) + B(\vtau,u) &= \int_{\partial_D} g_D\vtau\cdot \vn \, \d s, 
    && \text{for all }\vtau\in H_0(\text{div};D),\\[.5em]
  B(\vsigma,v) &= -\int_D f_\text{source}v\,\d x 
    && \text{for all  }v\in L^2(D)
\end{align*}
with the bilinear forms
\[
  A(\vsigma,\vtau) := \int_D a^{-1}\vsigma\cdot\vtau\,\d x
  \quad \text{ and } \quad
  B(\vsigma,u) := \int_D u(\nabla\cdot\vsigma)\,\d x.
\]
The (weak) differentiability properties required by this variational formulation lead us to seek the solution $(u,\vsigma)$ in the solution space $L^2(D) \times H_0(\text{div};D)$. 
We approximate the solution by introducing a triangular mesh ${\cal T}_h$ on $D$ and approximate $L^2(D)$ with the finite-dimensional subspace of piecewise constant functions on ${\cal T}_h$ and approximate  $H_0(\text{div};D)$ by the space of lowest-order Raviart-Thomas elements (cf.\ \cite{ErnGuermond2021}). 
The resulting system of discrete equations
\begin{equation*}
  \begin{bmatrix} \vA & \vB\\ \vB^\top & \vO \end{bmatrix}
  \begin{bmatrix} \vsigma^h \\u^h  \end{bmatrix}
  = 
  \begin{bmatrix} \vg_D^h \\ \vf_h \end{bmatrix}
\end{equation*}
has a sparse and symmetric system matrix and can be solved efficiently. 
For further details on mixed finite element methods see \cite{boffi2013mfem,ErnGuermond2021}.

Our quantity of interest is then approximated as
\begin{equation} \label{eq:QoI-flux}
    Q = Q(\vsigma(\omega)) \approx \int_{\partial_\text{out}} -\vsigma^h \cdot \vn \,\d s.
\end{equation}
The finite element discretization error for different types of quantities of interest for the stationary diffusion equation is analyzed in \cite{TeckentrupEtAl2013}.

The implementation of these finite element computations in \textsc{Matlab} is publicly available on \texttt{github}\footnote{\url{https://github.com/ToniKowalewitz/LearningToIntegrate}}. 
The three meshes employed were generated using \textsc{Matlab}'s \texttt{Partial Differential Equation Toolbox}.

\end{appendices}
\end{document}